\documentclass[review]{elsarticle}
\usepackage{lineno,hyperref}
\modulolinenumbers[5]
\usepackage{amsmath,amsfonts,amssymb,amsthm}
\usepackage[utf8]{inputenc}
\usepackage[english]{babel}
\usepackage{fullpage}
\usepackage{graphicx}
\usepackage{caption}
\usepackage{rotating}
\usepackage{longtable}
\usepackage{lscape}
\usepackage{pst-plot}
\usepackage{pst-all}
\usepackage{color, colortbl}
\usepackage{flexisym}
\usepackage{subcaption}
\usepackage{array}
\usepackage{tabularx}
\usepackage{multirow}
\usepackage{booktabs}
\usepackage{xcolor}
\usepackage{mwe}

\usepackage{enumitem}
\usepackage{multicol} 
\usepackage{tcolorbox} 
\usepackage{framed}
\usepackage{placeins}
\usepackage{amsmath,amssymb}

\usepackage{tikz}
\usepackage{pgfplots} 
\usepackage{tkz-fct}
\usepackage{mwe} 

\usepackage{pgf}
\usepackage[absolute, overlay]{textpos}
\textblockorigin{\paperwidth}{0.0 pt}

\pgfplotsset{compat=1.9}
\usepgfplotslibrary{statistics}
\usetikzlibrary{positioning}

\usepackage{setspace}
\newpsobject{malla}{psgrid}{subgriddiv=1,griddots=10,gridlabels=6pt}
\newcommand\numberthis{\addtocounter{equation}{1}\tag{\theequation}}
\makeatletter
\def\psls@dotteddashed{
[ 0 \psk@dotsep CLW add \psk@dash CLW add ] 0 setdash  1 setlinecap  stroke
}
\makeatother
\newcommand{\fixes}[1]{\color{black}{#1}\color{black}}

\begin{document}
\renewcommand{\tablename}{Table}

\journal{Journal of Sustainability Analytics and Modeling}
\begin{frontmatter}
\title{\fixes{On Carbon Taxes Effectiveness to Induce a Clean Technology Transition: An Evaluation Framework Based on Optimal Strategic Capacity Planning}} 
\author[1]{N. Wolf \corref{cor}}\ead{nathaliaiwolf@gmail.com}\cortext[cor]{Corresponding author.} 
\author[1]{P. Escalona}
\author[2]{A. Angulo}
\author[3]{J. Weston} 
\address[1]{Department of Industrial Engineering, Universidad T\'ecnica Federico Santa Mar\'ia, Avenida Espa\~na 1680, Valpara\'iso, Chile.}
\address[2]{Department of Electrical Engineering, Universidad T\'ecnica Federico Santa Mar\'ia, Avenida Espa\~na 1680, Valpara\'iso, Chile.}
\address[3]{Friedrich-Alexander-Universität Erlangen-Nürnberg, Department Mathematik, Cauerstraße 11, 91058 Erlangen, Germany.}
\begin{abstract}
\fixes{This paper studies carbon taxes effectiveness to induce a transition to cleaner production when a firm faces different technologies and demands. 
To determine carbon taxes effectiveness, we propose a framework based on a strategic capacity planning under carbon taxes model, that consider proper perfomance measures. 
The model, which is formulated as a mixed integer linear problem (MILP), considers issues that previous work have not studied jointly, such as machine replacement, workforce planning, and maintenance. 
The effectiveness measures consider levels of clean production and periods to reach a technological transition. 
Our computational experiments, based on a real case, have shown that carbon taxes by themselves do not necessarily induce a transition to clean production, since their effectiveness depends on the available technology relationship and the demand magnitude.  }
\end{abstract}
\end{frontmatter}
\allowdisplaybreaks
\section{Introduction}\label{sec:IN}

The balance of gases in the atmosphere creates what is known as the \textit{greenhouse effect}, which traps the sun's heat and, so far, has maintained the planet's temperature at suitable levels for the development of life.
However, in the last centuries, human activity has disrupted this balance \citep{cook2016consensus}.
Since the Industrial Revolution, thousands of tons of carbon emissions have been released into the atmosphere each year.
Carbon dioxide is one of the primary greenhouse gasses, and once its concentration in the atmosphere increases, more heat is trapped on the planet's surface, raising temperatures beyond normal levels.
The consequences of this disruption to the greenhouse effect, labeled as climate change, are already perceived through the ecosystem, and predictions state that there are less than two decades to reduce carbon emissions before reaching a no return point \citep{IPCC18}.
Therefore, the sustained increase in carbon emissions is of global concern \citep{yue2015high}.

In response to the threat of climate change, 195 countries signed the Paris Agreement in 2016, with the main goal of holding the increase in global temperatures below 2 °C pre-industrial levels.
However, efforts to accomplish this goal have proven to be insufficient, even more, it was estimated that by 2030 the global emissions will double what should be emitted to fulfill the controlled rise of temperatures \citep{UN2019}.

To control emissions, several policies and mechanisms have been designed to promote renewable energy and cleaner technologies \citep{zhang2021impact}.
Nowadays, carbon pricing is the main mechanism adopted by environmental authorities \citep{du2015expanding}, as they have proven to be an effective approach to reduce national emissions \citep{ghazouani2020exploring}. 
Two carbon pricing mechanisms can be distinguished. 
The first ones are \textit{carbon taxes}, which set a specific value per ton of carbon emissions, providing certainty in the cost faced by emitters. \fixes{Existing carbon taxes usually increase periodically, such as the French, Swedish, and Canadian carbon tax \citep{criqui2019carbon}. } 
The second ones are \textit{emissions trading systems} (ETS), which set a limit on the total tonnes of carbon emissions. ETS distribute the permitted emissions between emitters, who can later sell or buy to other emitters \citep{organizacion2013climate}.
Although both mechanisms promote a reduction of emissions \citep{he2015production}, 
carbon taxes are more swiftly implementable than ETS, making \fixes{them the preferred }option for environmental authorities \citep{dissanayake2020}. 

\fixes{The effectiveness of carbon taxes on national emissions, along with its effect on macro indicators has been widely studied, e.g. economic growth, employment and poverty. 
However, even when it is vital to reduce the environmental impact from the industrial sector \citep{huisingh2015recent}, research on the effectiveness of carbon taxes in individual firms is still scarce \citep{he2018upgrade}. 
We consider that an efficient mechanism to determine the effectiveness of carbon taxes on a firm is using the \textit{strategic capacity planning}. A strategic capacity planning under carbon taxes defines the optimal transition to clean technology of a firm, from which the effectiveness of carbon taxes can be measure in terms of the cleaning production level to reached, and the corresponding period of the transition. 
To the best of our knowledge, no work has considered a strategic capacity planning to study the optimal technological transition induced on a individual firm under carbon taxes. 

In this paper, we address the effectiveness of carbon taxes to induce a transition to cleaner production technologies on a firm that can choose between different technologies in terms of emissions and investment cost. 
The effectiveness of carbon taxes is determined using a framework based on a strategic capacity planning problem, which is formulated as a single-site and multi-item strategic capacity model under carbon taxes and deterministic demand. 
Our model combine capacity expansion decisions that previous work have not jointly consider, such as capacity expansion through investment and replacement, expansion or reduction through operational cost (workforce planning), demand satisfaction through production and inventory, allocation, maintenance, sale of discarded machines, and carbon taxes. 
Furthermore, we define measures for carbon taxes effectiveness, based on the transition levels and the periods in which clean production is reach. 
Thus, a carbon tax is effective if by the end of the planning horizon a technological transition is induced on the firm, the greater the transition to clean production and the sooner it occurs, the higher the effectiveness of the carbon tax.

The main contributions in this paper can be summarized as follows:
(i) We address the effectiveness of carbon taxes to induce a transition to cleaner technologies on an individual firm perspective; 
(ii) we proposed a framework to evaluate the effectiveness of the carbon tax based on a strategic capacity planning model;
(iii) we formulated and solve a strategic capacity planning model containing decisions not previously addressed jointly; and 
(iv) we determined factors that influence the effectiveness of carbon taxes.

The rest of the paper is structured as follows. 
In Section \ref{sec:RW}, a literature review on related work of carbon taxes effectiveness and strategic capacity planning is presented. 
In Section \ref{sec:fram} we present the proposed framework and formulate the strategic capacity planning under carbon taxes model. }
Section \ref{sec:CS} reports the computational experiments that explore carbon taxes effectiveness. Finally, Section \ref{sec:CO} presents our conclusions and future lines of research. 
\section{Literature review}\label{sec:RW}
\fixes{
There are two main streams of research related to this paper that need to be explored, (i) the effectiveness of carbon taxes, and (ii) the strategic capacity planning models.

The effectiveness of carbon taxes has been approached from different perspectives in literature. 
In this paper, we classified the effectiveness of carbon taxes under a \textit{macro}, \textit{supply chain}, and \textit{individual firm} perspective.  
The macro perspective focus on quantifying the effect that carbon taxes have on national or regional indicators \citep{hajek2019analysis, green2021does}, 
the supply chain perspective focus on the effect that carbon taxes have on the design of the logistic networks \citep{reddy2020effect}, and 
the individual firm perspective focus on how carbon taxes influence the operational decisions of firms \citep{zhou2020carbon}. 

From the macro perspective, \cite{jia2020rethinking}, and \cite{gokhale2021japan} analyze qualitatively the effectiveness of China and Japan's carbon taxation policy, respectively. Both authors determine that the value of the carbon tax, to decrease carbon emissions and transition out of fossil fuels, should be higher than actual values. 
\cite{lin2018energy}, and \cite{bernard2018effects} study the consequences that the carbon taxes have on the gross domestic product (GDP) of China and British Columbia, respectively. They conclude that carbon taxes do not cause a significant change on the GDP. 
\cite{renner2018poverty} and \cite{khastar2020does} simulated the possible impact of carbon taxes on social welfare and poverty in Mexico and Finland, showing a negative impact for high carbon taxes. 
\cite{brown2020all} and \cite{metcalf2020macroeconomic} study the implications of a carbon tax in employability in the U.S. and Europe. \cite{brown2020all} observed a positive relationship of the carbon tax and employment in the U.S., while \cite{metcalf2020macroeconomic} found no robust evidence of a negative impact on employment. 

For the supply chain perspective, a comprehensive review on the reduction of carbon emissions in supply chains can be found on \cite{das2018low} and \cite{chelly2019consideration}. 

In particular, \cite{zakeri2015carbon} and \cite{hariga2017integrated} research the impact of carbon taxes in the supply chain cost. 
\cite{zakeri2015carbon} compares the effect on costs of carbon taxes and ETS, observing an increase in costs under a carbon tax policy. 
\cite{hariga2017integrated} observes that incorporating carbon taxes in the supply chain, results in a minor increase in the operational costs. 
\cite{meng2018make} and \cite{cao2020production} study how carbon taxes affect production. 
\cite{meng2018make} observes that the cost for manufacturing firms is more sensitive to carbon taxes than firms that buy the product. 
\cite{cao2020production} concludes that production quantity decreases proportionally to the carbon tax value.
\cite{li2017production} and \cite{haddadsisakht2018closed} study the changes of transportation choices under carbon taxes. 
\cite{li2017production} observes that under increasing carbon taxes, the firm outsources more transportation services.  
\cite{haddadsisakht2018closed} observes that adjusting transportation capacities offers a better strategy than building extra facilities in response to carbon taxes uncertainty.
\cite{zhou2018pricing} and \cite{feng2020pricing} study the effect of carbon taxes on pricing, showing that final prices and wholesale prices increase as the carbon tax value increase. 
\cite{yi2018cost} and \cite{liu2021emission} explore the cooperation between the echelons of the supply chain. 
\cite{yi2018cost} observes that depending on the initial level of contamination and value, a carbon tax policy can promote cooperation in the supply chain. 
\cite{liu2021emission} results shows that,when manufacturers and retailers cooperate under carbon tax regulations, the emissions reduction and profit increment can be optimized.

Regarding the effectiveness of carbon taxes to reduce emissions in individual firms, a comprehensive review on carbon emissions in operational decisions can be found on \cite{zhou2020carbon}. 
\cite{lin2017pull} and \cite{tang2018reducing} study the impact of carbon taxes on inventory decisions.  
\cite{lin2017pull} observes that the number of shipments and the lot size depend on the carbon tax policy, a progressive carbon tax policy offers flexibility in deciding the number of shipments and the sizes, and a multiple delivery policy is better than a single delivery policy. 
\cite{tang2018reducing} shows that emissions from holding and transportation can be reduced changing the lot size and reorder point. 
\cite{ma2018optimal} and \cite{kondo2019green} study how supplier selection change under carbon taxes. 
Both authors observe that carbon taxes can induce manufacturers to choose lower emitter supplier once carbon tax value is increase enough. 
\cite{wang2018optimization}, \cite{qin2019optimization} and \cite{li2020investigating} study the effect of carbon taxes on the optimal path of a vehicle routing problem.  
Their results show that a carbon tax policy, under a minimization costs approach, can reduce carbon dioxide emissions in logistics network. 
Furthermore, the effect of the speed and capacity of the vehicles becomes more relevant as the carbon taxes increase. 
To the best of our knowledge, only \cite{drake2016technology} and \cite{song2017capacity} study the impact of carbon taxes on the firm technology and capacity expansion from an individual firm perspective. However, they use a model based on Stackelber games that only allows them to analyze the capacity expansion in the short term. 
In summary, no previous work explores the long term effectiveness of carbon taxes to induce a transition toward cleaner technologies in a firm following a strategic planning perspective. 
}

\fixes{The second literature stream related to this paper, is strategic capacity planning. }
A comprehensive review can be found in \citet{verter1992integrated}, \cite{van2003commissioned}, \cite{wu2005managing}, \cite{julka2007review}, and the most recent in \cite{martinez2014review}. 
\cite{martinez2014review} classified the strategic capacity models according to the number of sites involved in the process (single-site or multi-site), the type of capacity considered (investing, outsourcing, reduction, and replacement) and the uncertainty considered in the formulation (deterministic or stochastic).
Furthermore, they described major decisions that a strategic capacity planning problem should address: capacity size, capacity location, allocation, technology selection, production and inventory, backorders, workforce planning, new product development, and financial planning. 

The capacity expansion problem of \citet{drake2016technology}, studies a single-site and single-item capacity model under uncertain demand.  
They considered the capacity expansion through investment by acquiring units of capacity of two different technologies (dirty or clean) and demand is satisfied through production.
The model is formulated as a two-stage stochastic linear programming which maximizes the expected profit.
They derived an optimal closed-form solution that depends on the problem parameters.
\citet{song2017capacity}, extend the work of \citet{drake2016technology}, considering capacity  at the beginning of the planning horizon. 
However, the two related works of \citet{drake2016technology} and \citet{song2017capacity} left out crucial decisions for the strategic capacity plan with carbon taxes, namely, machine replacement and workforce planning. 

Carbon taxes promote a transition to clean technologies. 
Therefore, it is essential to consider discarding old machinery when planning capacity under carbon taxes.
A literature review of equipment replacement is presented by \cite{hartman2014equipment}.
\cite{chand2000model} analyzed a single-site and single-item strategic capacity model under deterministic demand. 
They considered capacity expansion through investment and replacement of machines, where maximum periods of use is allowed before a machine must be replaced.  
The model is formulated as a MILP and solved by a heuristic procedure.
\cite{mitra2014optimal} studied a single-site and multi-item strategic capacity model under uncertain demand.
They considered capacity expansion through investment and replacement of machines, where an upgrade of equipment constitutes a machine replacement, and demand is satisfied by production and inventory. The model is formulated as a two-stage stochastic MILP and solved using a commercial solver.
\cite{benedito2016single} formulated a single-site and multi-item strategic capacity under deterministic demand. They considered expansion through investment and replacement of machines, demand is satisfied through production and inventory, and the selling of machines associated with the salvage value of the equipment. The model is formulated as a MILP and solved using a commercial solver. 
\cite{wang2017capacity} analyzed a single-site and multi-item strategic capacity under deterministic demand.
They considered capacity expansion through investment and replacement of machines, with technology selection, where demand is satisfied through production and inventory. 
The model is formulated as a MILP, using a stochastic dynamic approach to address technology uncertainty, and solved by a genetic algorithm.

\fixes{Regarding the workforce planning, for capital intensive firms, the implementation of extra work shifts is an additional tool to increment capacity and delay the costs of machine acquisition. Few papers consider the workforce in strategic capacity planning. }
\cite{fleischmann2006strategic} studied a multi-site and multi-item strategic capacity model under deterministic demand.
They considered capacity expansion through investment and replacement of machines, and through operational cost with overtime and change of work shifts. Additionally, they considered production planning with allocation of machines to satisfy demand and included the location as a decision variable.
The model is formulated as a MILP and solved using a commercial solver.
\cite{bihlmaier2009modeling} analyzed a multi-site and multi-item strategic capacity model under uncertain demand.
They considered expansion through investment by acquiring machines, and \fixes{through } operational costs by \fixes{adding } new work shifts, such as a late, night, \fixes{and }Saturday work shift.   
Additionally, they considered production planning with allocation and backorders to satisfy demand and included location as a decision variable.
The model is presented as a two-stage stochastic MILP and solved using a commercial solver. 
\cite{weston2017strategic} studied  single-site and multi-item under uncertain demand.
They considered capacity expansion through investment and operational cost by deciding the number of work shifts, and considered production planning with allocation to satisfy demand. 
Their model is formulated as a MILP, which presents a robust approach to include the uncertainty and was solved using a commercial solver. 

\fixes{
All the listed works, related to strategic capacity planning, illustrate the different characteristics considered. These characteristics can be grouped into (i) problem setting, (ii) decisions included, and (iii) formulation. For the problem setting we consider sites (single-site or multi-site), items (single-item or multi-item) and demand (deterministic or uncertain). In the decision we include the presented by \citet{martinez2014review}, and add three more, namely, maintenance of machines, sale of discarded equipment, and carbon taxes. With respect to formulation, we focus on capacity (discrete or continuous expansion), type of model, and solution approach. 
Table \ref{tabla_lr} shows the characteristics of the models presented in this section, as well as the model proposed in this work. 

\renewcommand{\thefootnote}{\alph{footnote}}
\setcounter{footnote}{0}

\begin{table}[h!]
\begin{minipage}[c]{\textwidth}
  \centering
  \caption{Characteristics of related strategic capacity planning models.} \label{tabla_lr} 
    \resizebox{\textwidth}{!}{    \begin{tabular}{l|ccc|crccccrcrrccc|ccc}
    \toprule
          & \multicolumn{3}{c|}{Setting} & \multicolumn{13}{c|}{Decisions}                                                                       & \multicolumn{3}{c}{Formulation} \\
    \midrule
    Authors and year & \begin{sideways}Single-site/Mult-site\footnotemark \end{sideways} & \begin{sideways}Single-item/Multi-item\footnotemark\end{sideways} & \begin{sideways}Demand\footnotemark\end{sideways} & \begin{sideways}Capacity size\footnotemark\end{sideways} & \multicolumn{1}{c}{\begin{sideways}Capacity Location\end{sideways}} & \begin{sideways} Allocation\end{sideways} & \begin{sideways}Technology Selection\end{sideways} & \begin{sideways}Production planning\end{sideways} & \begin{sideways}Inventory\end{sideways} & \multicolumn{1}{c}{\begin{sideways}Backordes\end{sideways}} & \begin{sideways}Workforce Planning\end{sideways} & \multicolumn{1}{c}{\begin{sideways} New Product Development\end{sideways}} & \multicolumn{1}{c}{\begin{sideways}Financial Planning\end{sideways}} & \begin{sideways}Maintenance\end{sideways} & \begin{sideways}Sale of discarded \end{sideways} & \begin{sideways}Carbon tax\end{sideways} & \begin{sideways} Capacity\footnotemark \end{sideways} & \begin{sideways}Type of model \end{sideways} & \begin{sideways}Solution Approach\end{sideways} \\
    \midrule
    \citet{chand2000model} & S     & S     & D     & E/R   &       &       &       &       &       &       &       &       &       &       &       &       & D    & MILP  & Heuristic \\
    \citet{fleischmann2006strategic} & M     & M     & D     & E/R   & \multicolumn{1}{c}{\checkmark} & \checkmark    &       & \checkmark    &       &       & \checkmark    &       &       &       &       &       &   C    & MILP  & Solver \\
    \citet{bihlmaier2009modeling} & M     & M     & U     & E/R   & \multicolumn{1}{c}{\checkmark} & \checkmark    &       & \checkmark    &       & \multicolumn{1}{c}{\checkmark} & \checkmark    &       &       &       &       &       &   C    & MILP  & Solver \\
    \citet{escalona2012expansion}  & S     & M     & D     & E     &       & \checkmark    &       & \checkmark    &       &       & \checkmark    &       &       &       &       &       & D    & MINLP & Solver \\
    \citet{mitra2014optimal} & S     & M     & U     & E/R   &       &       &       & \checkmark    & \checkmark    &       &       &       &       &       &       &       & D    & MILP  & Solver \\
    \citet{drake2016technology} & S     & S     & U     & E     &       &       & \checkmark    & \checkmark    &       &       &       &       &       &       &       & \checkmark    &    C   & LP    & Closed-form \\
    \citet{benedito2016single} & S     & M     & D     & E/R   &       &       &       & \checkmark    & \checkmark    &       &       &       &       & \checkmark    & \checkmark    &       &    C   & MILP  & Solver \\
    \citet{song2017capacity} & S     & S     & U     & E     &       &       & \checkmark    & \checkmark    &       &       &       &       &       &       &       & \checkmark    &   C   & LP    & Closed-form \\
    \citet{wang2017capacity} & S     & M     & D     & E/R   &       & \checkmark    & \checkmark    & \checkmark    & \checkmark    &       &       &       &       &       &       &       &   C     & MILP  & Heuristic \\
    \citet{weston2017strategic} & S     & M     & U     & E     &       & \checkmark    &       & \checkmark    &       &       & \checkmark    &       &       &       &       &       &  D    & MILP  & Solver \\
    \textcolor[rgb]{ 1,  0,  0}{Proposal} & \textcolor[rgb]{ 1,  0,  0}{S} & \textcolor[rgb]{ 1,  0,  0}{M} & \textcolor[rgb]{ 1,  0,  0}{D} & \textcolor[rgb]{ 1,  0,  0}{E/R} & \textcolor[rgb]{ 1,  0,  0}{} & \textcolor[rgb]{ 1,  0,  0}{\checkmark} & \textcolor[rgb]{ 1,  0,  0}{\checkmark} & \textcolor[rgb]{ 1,  0,  0}{\checkmark} & \textcolor[rgb]{ 1,  0,  0}{\checkmark} & \textcolor[rgb]{ 1,  0,  0}{} & \textcolor[rgb]{ 1,  0,  0}{\checkmark} & \textcolor[rgb]{ 1,  0,  0}{} & \textcolor[rgb]{ 1,  0,  0}{} & \textcolor[rgb]{ 1,  0,  0}{\checkmark} & \textcolor[rgb]{ 1,  0,  0}{\checkmark} & \textcolor[rgb]{ 1,  0,  0}{\checkmark} & \textcolor[rgb]{ 1,  0,  0}{D} & \textcolor[rgb]{ 1,  0,  0}{MILP} & \textcolor[rgb]{ 1,  0,  0}{Solver} \\
    \bottomrule
    \end{tabular}%
    }
    \footnotetext[1]{S (Single-site), M (Multi-site)}
    \footnotetext[2]{S (Single-item), M (Multi-item)}
    \footnotetext[3]{D (Deterministic), U (Uncertain)}
    \footnotetext[4]{E (Expansion), R (Replacement)}
    \footnotetext[5]{C (Continuous), D (Discrete)}
\end{minipage}
\end{table}

From Table \ref{tabla_lr} it can be seen that the strategic capacity planning under carbon taxes is still narrow, since current works with carbon taxes have not yet considered relevant decisions to determine the optimal transition to cleaner technologies, such as replacement, workforce planning, inventory, allocation or maintenance. 
All these decisions are addressed in this work, where we propose an integrated model which considers relevant decisions for a strategic capacity plan to evaluate the effect of carbon taxes on the clean technology transition of a firm. 
}

\fixes{
\section{Framework to determine the effectiveness of carbon taxes in a firm} \label{sec:fram} 
This section presents a framework to determine the effectiveness of carbon taxes to induce a transition to clean technology on a manufacturing firm.  
The framework is based on a strategic capacity planning model that allows evaluating different demand behavior, planning horizons, and types of technologies. The model considers relevant decisions on the strategic capacity planning literature, which have not been considered when evaluating the effect of carbon taxes on a firm. 
Furthermore, we propose two measures of effectiveness, the \textit{transition level} and \textit{transition period}. 
The structure of the proposed framework is illustrated in Figure \ref{fig:frame}, and will be further explain in what follows.

\begin{figure}[h!]
    \centering
    \includegraphics[scale=1.05]{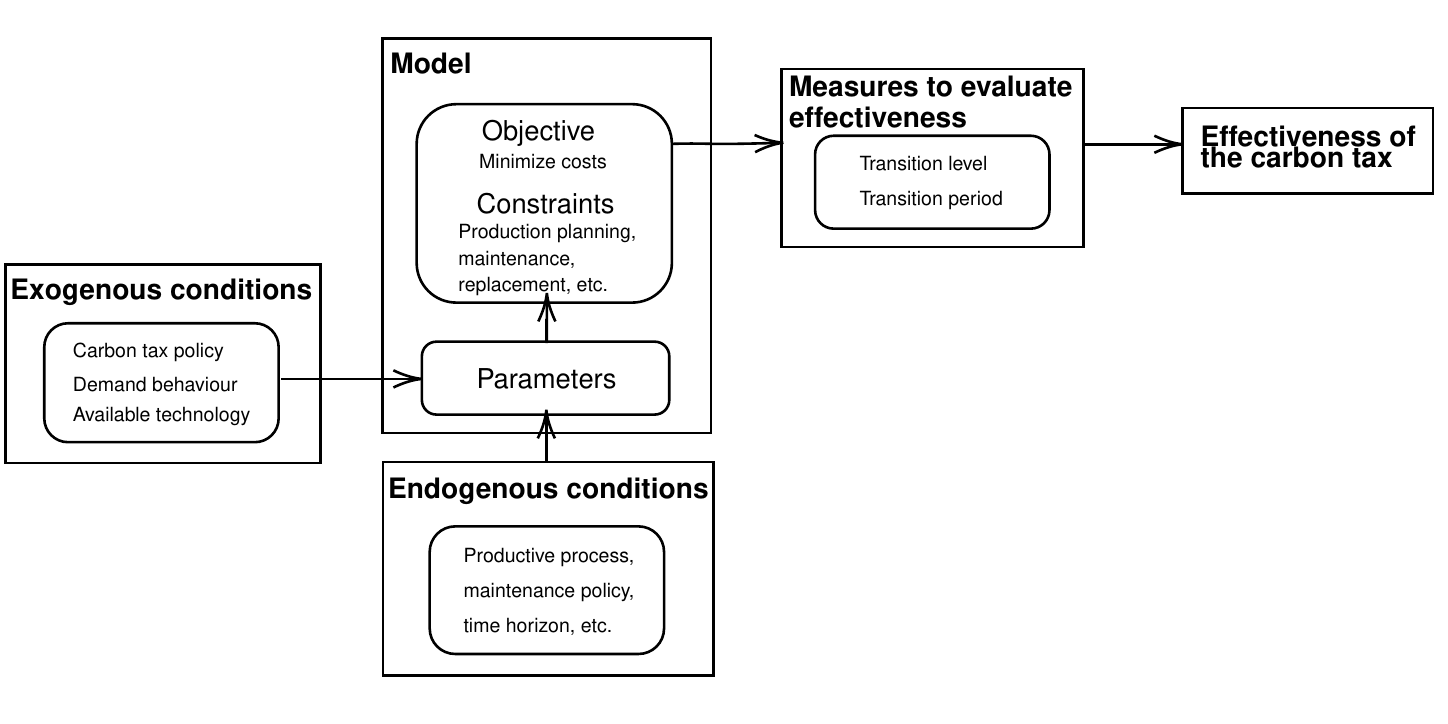}
    \caption{Structure of the framework to determine the effectiveness of carbon taxes in a firm.}
    \label{fig:frame}
\end{figure}

\subsection{Endogenous and exogenous conditions}

This framework allow us to evaluate carbon taxes effectiveness to induce a clean technology transition on a manufacturing firm, under different endogenous and exogenous conditions. 

The endogenous parameters of the stratetegic capacity planning are those which are specific to a firm. Different firms and production processes influence several operational costs, such as production, maintenance, and labor.  
Thus, the firm and the scope of its manufacturing process defines the endogenous parameters of the problem. 

On the other hand, the exogenous parameters depend on the external scenario in which the firm is immersed. 
In this work, we consider three factors that define the scenario faced by the firm. 
The first factor is the carbon tax policy, i.e., the magnitud and update frequency of the carbon tax implemented. 
The second factor is the demand faced by the firm, since the rise of demand is directly tied with carbon emission growth \citep{drake2013om} and the strategic capacity plan of a firm. 
The third factor is the technology available in the market, which differs in terms of the emissions by unit produced and the investment cost. 
}
\subsection{Strategic capacity planning model}\label{sec:MF}
Consider the strategic capacity planning of a process over a horizon of $T\ (t=1,...,|T|)$ periods, where increasing carbon taxes are imposed by the environmental authority.
In this process,  $I\ (i=1,...,|I|)$ items are produced, where an item correspond to a component or final product.
Let $d_{it}$ be the demand of item $i$ satisfied by production or inventory in period $t$. 
At each period, the planning of the process capacity can be expanded through the acquisition of machines (expansion through investment) or by modifying the number of work shifts (expansion through operational cost). 

Let $K^0$ be the set of operating machines at the beginning of the planning horizon ($t=0$) and $K^c$ be the set of machines that can be acquired over time. Thus, $K=K^0 \cup K^c \ (k=1,...,|K|)$ is the set of machines, where $K^0$ can be an empty set.
In this paper, machine refers to equipment as well as production line. The machine capacity is based on available time by period, which depends on the operational conditions of the process and the number of work shifts at the period. Machines wear out on use and are maintained under a preventive policy. Without loss of generality, it is considered that the policy is \textit{as good as new}, and is executed on working hours.
Furthermore, as carbon taxes \fixes{are expected to } induce the replacement of older machines with cleaner ones, the discard (and sale) of machines is considered. 

Since the number of work shifts can be modified in each period, by opening or closing work shifts, it is possible to hire or fire worker, according to labor requirements.  
Moreover, at each period the process can use one, two, or three work shifts with the same length. 

The objective of our model is to determine the strategic capacity plan, \fixes{that } minimizes costs over the planning horizon considering increasing carbon taxes. In this paper, the strategic capacity plan costs are investment, production, maintenance, labor, hiring-firing, change of work shifts, inventory, and carbon taxes; as well as the revenue from discarded (sold) machines. 

\subsubsection{Machine state diagram}
The operational state of a machine and its change of state between periods (transitions) can be efficiently represented with a state diagram. 
Let $S = \{0,1,2,3\}$ be the set of \textit{operational states} of a machine at any period of the planning horizon, where 
$s=0$ means that the machine is not operating since it has not been purchased or has already been discarded;
$s=1$ means that the machine operates one work shift;
$s=2$ means that the machine operates two work shifts; and
$s=3$ means that the machine operates three work shifts. 
Let $E = \{(s,s') \in S \times S \}$ be the set of \textit{state transitions} that can occur at the beginning of any period of the planning horizon.
Thus, the state of the machine during a period $t \in T$ is defined by the state transition $e=(s,s') \in E$, where $s$ is the state during period $t-1$, and $s'$ is the state during period $t$.
We label all elements of the state transition set as follows: $E=\{e_0, ... , e_{15}\}$, where $e_0=(0,0)$, $e_1=(0,1)$, $e_2=(0,2)$, $e_3=(0,3)$, $e_4=(1,0)$, $e_5=(1,1)$, $e_6=(1,2)$, $e_7=(1,3)$, $e_8=(2,0)$, $e_{9}=(2,1)$, $e_{10}=(2,2)$, $e_{11}=(2,3)$, $e_{12}=(3,0)$, $e_{13}=(3,1)$, $e_{14}=(3,2)$ and $e_{15}=(3,3)$. 

\begin{figure}[!htb]
{
\begin{subfigure}[t]{0.23\textwidth}
    \centering
    \includegraphics[width=0.3\linewidth]{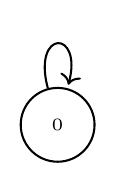} 
    \caption{\scriptsize $E^0$} 
    \label{fig:e0} 
    \vspace{4ex}
  \end{subfigure} 
  \begin{subfigure}[t]{0.23\textwidth}
    \centering
    \includegraphics[width=0.6\linewidth]{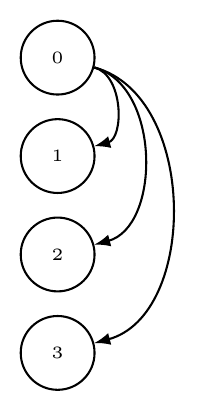}  
    \caption{\scriptsize $E^1$} 
    \label{fig:e1} 
    \vspace{4ex}
  \end{subfigure} 
  \begin{subfigure}[t]{0.23\textwidth}
    \centering
    \includegraphics[width=0.85\linewidth]{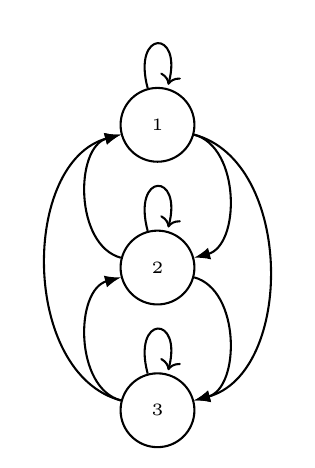} 
    \caption{\scriptsize $E^2$} 
    \label{fig:e2} 
  \end{subfigure} 
  \begin{subfigure}[t]{0.23\textwidth}
    \centering
    \includegraphics[width=0.6\linewidth]{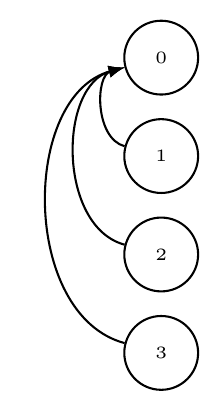}  
    \caption{\scriptsize$E^3$} 
    \label{fig:e3} 
  \end{subfigure} 
  
}
\caption{State diagram partitions.} \label{fig:E}
\end{figure}

The set of state transitions $E$ can be conveniently partitioned into four subsets.
Let $E^0 = \{e_{0} \}$ be the transition \fixes{set }that represents a machine that is \fixes{currently } inoperative, since it has not been purchased or has already been discarded;
$E^1 = \{e_{1},$ $e_{2}$, $e_{3} \}$ be the transitions \fixes{set }that represents the entry into operation of a machine that has not been purchased; 
$E^2=\{e_5$, $e_6,$ $e_7,$ $e_{9},$ $e_{10},$ $e_{11},$ $e_{13},$ $e_{14},$ $e_{15} \}$ be
the transition \fixes{set } that represents a machine that was operating the previous period and \fixes{stays } in a state different from zero; and 
$E^3 = \{e_4$, $e_8,$ $e_{12}\}$ be the transition \fixes{set } that represents the discard of a machine that was operating during the previous period.
The state diagrams representing the four partitions of $E$ are illustrated in Figure \ref{fig:E}.

Figure \ref{fig:e0} shows a machine that during period $t-1$ is in state $s=0$ and in period $t$ remains in state $s'=0$.
Figure \ref{fig:e1} shows a machine that during period $t-1$ is in state $s=0$ and in period $t$ is in a state $s' \neq 0$, i.e., the machine begins to operate at the beginning of period $t$.
Figure \ref{fig:e2} shows a machine that during period $t-1$ is in a state $s \neq 0$ and in period $t$ is also in a state $s' \neq 0$, i.e., machine stays in the same state or change to another state different from $0$.
Figure \ref{fig:e3} shows a machine that during $t-1$ is in a state $s \neq 0$ and in period $t$ is in state $s' = 0$, i.e., the machine is discarded at the beginning of period $t$. 

\subsubsection{The strategic capacity problem with carbon taxes}

Under a state approach, a machine remains or changes its state at the beginning of a period. Let $X_{ke}^{t}$ be equal to 1 if \fixes{state transition }$e \in E$ occurs at the beginning of period $t$ for the machine $k \in K$; 0 otherwise. 
Furthermore, for any $e=(s,$ $s') \in E$ in period $t$ we define $H(e)=s'$ as the state of machine $k$ in period $t$, and $T(e)=s$ as the state of machine $k$ in period $t-1$. 

To establish the relationship between the transitions of any machine throughout the planning horizon, the following set of constraints are defined: 
\begin{align}
\sum_{e\in E: T(e)=s}{X_{ke}^{t+1}} - \sum_{e\in E: H(e)=s}{X_{ke}^t} &= 0 && \forall k\in K, s \in S, t \in T \setminus \{|T|\}, \label{1.1} \\
\sum_{e \in E}{X_{ke}^t} &= 1  && \forall k\in K, t \;\fixes{\in T \setminus \{1\}}, \label{1.2} \\
\sum_{t \in T}\sum_{e \in E^1} X_{ke}^t &\leq 1 && \forall k \in K^c, \label{1.3} \\
X_{ke}^t &\in \{0,1\} && \forall k \in K, e \in E, t \in T. \label{1.4}
\end{align}

Constraints (\ref{1.1}) ensures the continuity of transitions for any machine throughout the planning horizon.  
Constraints (\ref{1.2}) \fixes{forces } that for any machine only one transition occurs at the beginning of any period. 
Constraints (\ref{1.3}) ensures that during the planning
horizon, a machine $k \in K^c$ can only be bought once, i.e., for any machine $k \in K^c$ a transition $e \in E^1$ can occur, at most, once. 

It is considered that operating machines are used the same number of work shifts over the planning horizon, i.e., at the beginning of period $t$ \fixes{all the }operating machines transition to the same state $s'$. 
We assumed that the state $s$ of every machine $k \in K^0$ at the beginning of the planning horizon ($t=0$) is a known parameter, denoted by $s_0$. In the same way, the initial state of a machine $k \in K^{c}$ is $s=0$.
Let $Z_{s}^t$ be equal to 1 if $s$ work shifts are used during period $t$; 0 otherwise. 
To establish the initial state of the different machines and their relation to the work shifts used, the following constraints are defined:
\begin{align}
\sum_{e\in E^2 \cup E^3: T(e)=s_0}{X_{ke}^1} &= 1 && \forall k\in K^{0},  \label{1.5} \\
\sum_{e\in E^0 \cup E^1}{X_{ke}^1} &= 1 && \forall k\in K^{c}, \label{1.6} 
\\
\sum_{s\in S}{Z_{s}^t} &= 1  && \forall t \in T, \label{1.7}
\\
\sum_{e\in E^1 \cup E^2: H(e)=s}{X_{ke}^t} &\leq Z_{s}^t && \forall s \in S, k\in K, t \in T, \label{1.8} \\
Z_{s}^{t} &\in \{0,1\} && s \in S, t \in T. \label{1.9}
\end{align}

Constraints (\ref{1.5}) ensures that for any machine $k \in K^0$ only one transition $e \in E^2 \cup E^3$ can occur at the beginning of period $t=1$, since the state of this kind of machines in $t=0$ is $s=s_{0}$. 
Constraints (\ref{1.6}) ensures that for any machine $k \in K^c$ only one transition $e \in E^0 \cup E^1$ can occur at the beginning of period $t=1$, since the state of this kind of machines in $t=0$ is $s=0$.
Constraints (\ref{1.7}) and (\ref{1.8}) follows that all operating machines \fixes{are used for the same number of work shifts in period $t$}.

The opening and closing of a work shift at the beginning of period $t$ is defined by the binary variable $O_{st}$ and $C_{st}$, respectively. Let $O_{st}$ ($C_{st}$) be equal to $1$ if work shift $s$ is opened (closed) at the beginning of period $t$; $0$ otherwise. Thereby, the opening and closing of a work shift is defined as follows, 
\begin{align}
Z_{s}^t-Z_{s}^{t-1} &\leq O_{s}^t && \forall s \in S, t \in T, \label{1.10} \\ 
Z_{s}^1-Z_{s_0}^{0} &\leq O_{s}^1 && \forall s \in S, \label{1.33} \\
Z_{s}^{t-1}-Z_{s}^{t} &\leq C_{s}^{t}  && \forall s \in S, t \in T, \label{1.11} \\ 
Z_{s_0}^{0}-Z_{s}^{1} &\leq C_{s}^{1}  && \forall s \in S, \label{1.34} \\
O_{s}^t, C_{s}^{t} &\in \{0,1\} &&\forall s \in S, t \in T, \label{1.12}
\end{align}
where $Z_{s_0}^{0}$ is a known parameter that is equal to 1 if $s_0$ work shift are used in period $t=0$; 0 otherwise. 

Constraints (\ref{1.10}) to (\ref{1.34}) ensures that if there is a change of work shift between $t-1$ and $t$, an opening and closing of work shift will occur.

As mentioned previously, the demand for each period is satisfied with production and/or inventory.
Let $Y_{ik}^{t}$ be the quantity of item $i$ produced by machine $k$ in period $t$, and $I_{i}^{t}$ be the on-hand inventory of item $i$ at the end of period $t$. Thus, the following constraints ensure that the demand is met: 
\begin{align}
I_{i}^{t-1} + \sum_{k\in K}{Y_{ik}^t} - I_{i}^t &\geq d_{it} &&  \forall i\in I, t \; \fixes{\in T \setminus \{1\}}, \label{1.13} \\
I_{i}^{0} + \sum_{k\in K}{Y_{ik}^1} - I_{i}^1 &\geq d_{i1} &&  \forall i\in I, \label{1.32} \\ 
Y_{ik}^{t} &\geq 0 &&\forall i \in I, k \in K, t \in T, \label{1.14} \\
I_{i}^{t} &\geq 0 &&\forall i \in I, t \in T, \label{1.15}
\end{align}
where $I_{i}^{0}$ is the on-hand inventory of item $i$ at the end of period $t=0$, which may be zero. 
It should be noted that other sources of production used in production planning, such as backorders, outsourcing, or extra time are not considered. 

To consider the wear of machine due to their use, we assumed, without loss of generality, a fixed-time maintenance policy, where $FTM_k$ is the maximum fixed time between maintenance of machine $k$. 
Let $W$ $(w=1,...,$ $|W|)$ be the numbers of maintenance that can be carried out for each machine, where the time required for the $w$\textit{-th} maintenance of machine $k$ is known and denoted by $RMT_{wk}$. Let $M_{wk}^t$ be equal to 1 if the $w$\textit{-th} maintenance of the machine $k$ is performed in period $t$; 0 otherwise, and $TM_{k}^t$ be the accumulated production time of machine $k$ from its last maintenance to the beginning of period $t$.
Consequently, the following sets of constraints define the execution of maintenance:
\begin{align}
TM_k^t &\geq TM_{k}^{t-1} + \sum_{i\in I}{\frac{Y_{ik}^{t-1}}{r_{ik}}} - FTM_{k} \sum_{w \in W}{M_{wk}^{t}} && \forall k\in K, t \in T, \label{1.16}
\\ 
TM_{k}^{t} &\leq FTM_{k} && \forall k\in K, t \in T, \label{1.17}
\\
\sum_{\tau=1}^{t}{M_{wk}^{\tau}} &\geq M_{w+1k}^{t} && \forall k\in K, t \in T, w \in 1, ... \lvert W \rvert -1, \label{1.18}
\\
\sum_{t\in T}{M_{wk}^{t}} &\leq 1  && \forall k\in K, w \in W, \label{1.19} 
\\
{M_{wk}^{t}} &\leq \sum_{\tau = 1}^{t}\sum_{e \in E^1 \cup E^2}{X_{ke}^{\tau}}  && \forall k\in K, t \in T, w \in W, \label{1.20} \\
TM_{k}^{t} &\geq 0 &&\forall k \in K, t \in T, \label{1.21} \\
M_{wk}^{t} &\in \{0,1\} &&\forall w \in W, k \in K, t \in T, \label{1.22}
\end{align}
where $r_{ik}$ is the production rate of item $i$ produced with machine $k$. 

Constraints (\ref{1.16}) updates the time counter from the last maintenance and constraints (\ref{1.17}) ensures that this time does not exceed the maximum fixed time \fixes{between } maintenance of machine $k$. Constraints (\ref{1.18}) implies that the $w$\textit{+1}\textit{-th} maintenance cannot be executed until all previous maintenance have occurred. Constraints (\ref{1.19}) ensure that the $w$\textit{-th} maintenance can be executed, at most, once. Constraints (\ref{1.20}) ensures that a maintenance is only applied to operating machines, i.e., only machines that at the beginning of period $t$ stay or transition to a state $s \neq 0$ can be maintained. 

Machines have a useful life which wears out according to its production time. 
Let $RL_k^{t}$ be the remaining useful life of machine $k$ at the beginning of period $t$. 
Thus, the production time of a machine cannot exceed its remaining useful life.  
On the other hand, it is clear that the production and maintenance time cannot exceed the available time of the period. 
Consequently:  
\begin{align}
\sum_{i\in I}{\frac{Y_{ik}^{t}}{r_{ik}}} + \sum_{w \in W}{RMT_{wk}M_{wk}^{t}} &\leq \mu_{k} \sum_{s \in S}\sum\limits_{\substack{\; e \in E^1 \cup E^2:\\ H(e)=s}}  ls X_{ke}^{t}  && \forall k\in K, t\in T, \label{1.23}\\
\sum_{i\in I}{\frac{Y_{ik}^t}{r_{ik}}} &\leq RL_{k}^{t} && \forall k\in K, t\in T, \label{1.24}\\
RL_{k}^{t} &= RL_{k}^{t-1} + v_{k} \sum_{e\in E^1}{X_{ke}^{t}} - \sum_{i\in I}{\frac{Y_{ik}^{t-1}}{r_{ik}}} && \forall k\in K, t \;\fixes{\in T \setminus \{1\}}, \label{1.25} \\
RL_{k}^{1} &= RL_{k}^{0} + v_{k} \sum_{e\in E^1}{X_{ke}^{1}} && \forall k\in K, \label{1.31} \\
RL_{k}^{t} &\geq 0 &&\forall k \in K, t \in T,\label{1.26} 
\end{align}
where $\mu_{k}$ is the maximum utilization of machine $k$, $v_k$ is the useful life of machine $k$, and $l$ is the work shift length. Also, $RL_{k}^{0}$ represents the remaining useful life of machine $k$ at the end of period $0$, with $RL_{k}^{0}=0$ for any $k \in K^c$. 

Constraints (\ref{1.23}) ensures that for any machine the production and maintenance time is less or equal that its available time at any period. Constraints (\ref{1.24}) ensures that production time of a machine does not exceed the remaining useful life. Constraints (\ref{1.25}) and (\ref{1.31}) updates and initializes the remaining useful life of a machine, respectively.

As stated in Section \ref{sec:RW}, the replacement of machines is an essential decision to consider in a strategic capacity plan under carbon taxes.
The resale revenue depends on variable $RF_{k}^{t}$, defined as the remaining useful life of machine $k$ when discarded at the beginning of the period $t$. Consequently, the following constrains are defined: 
\begin{align}
RF_{k}^{t} &\leq v_{k}\sum_{e\in E^3}{X_{ke}^{t}} && \forall k \in K, t \in T, \label{1.27}
\\
RF_{k}^{t} &\leq RL_{k}^{t} && \forall k \in K, t \in T, \label{1.28} \\
RF_{k}^{t} &\geq 0 &&\forall k \in K, t \in T. \label{1.29}
\end{align}

Constraints (\ref{1.27}) and (\ref{1.28}) \fixes{condition } the value of the remaining life of a machine when discarded at the beginning of the period $t$.

The objective function involves minimizing the total cost of the strategic capacity plan, including the carbon taxes. 
We consider eight costs associated with the strategic capacity planning. \fixes{Morever}, since the discard of machines is allowed, the objective function also includes the revenue from selling the discarded machines. 
Consequently, we have:
$(i)$ the investment cost of acquiring machine $k$ in period $t$ ($CI_{k}^{t}$);
$(ii)$ the production cost for an item $i$  produced by machine $k$ in period $t$ ($CP_{ik}^{t}$);
$(iii)$ the cost of the $w$\textit{-th} maintenance of machine $k$ in period $t$ ($CM_{wk}^{t}$); 
$(iv)$ the unitary labor cost for machine $k$ in period $t$ ($CL_k^{t}$);
$(v)$ the cost of hiring ($CA_{k}^{t}$) and firing ($CF_{k}^{t}$) workers;
$(vi)$ the cost of opening ($CO_{s}^{t}$) and closing ($CC_{s}^{t}$) work shifts $s$;
$(vii)$ the holding cost per unit and unit period of item $i$ in period $t$ ($CH_{i}^{t}$);
and, $(viii)$ the value of carbon taxes in period $t$ ($CT_{t}$), which increases periodically, according to the carbon tax policy adopted by the environmental authority. 
Furthermore, we consider the revenue from sold machines. 
Let $\alpha_{k}^{t}$ be the price per hour of useful life remaining of machine $k$ that was discarded in period $t$. 
All the aforementioned costs are properly brought to their current value.

Carbon taxes take various forms depending on the regulatory policies of emission control. Similarly to \cite{benjaafar2013} and \cite{mohammed2017multi}, we consider that carbon emissions come from the energy and the resources used in production and inventory.
Let $ep_{ik}$ be the emissions for producing one unit of item $i$ using machine $k$, 
and let $eh_{i}$ be the emissions of holding per unit and unit period of item $i$.
\fixes{Therefore}, the carbon emissions in period $t$ are $\mathcal{E}^t = \sum_{i \in I}(\sum_{k \in K} {ep_{ik}  Y_{ik}^{t}} + eh_{i} I_{i}^{t})$. Thus, the strategic capacity planning under carbon taxes \fixes{(\textbf{SPT}) } is defined as the following MILP: 
\begin{align*}
 \textbf{SPT}:
 \quad \min & \;\; \sum_{k\in K, t\in T} \left\lbrace \sum_{e\in E^1}{CI_{k}^{t} X_{ke}^{t}} + \sum_{i\in I}{CP_{ik}^{t} Y_{ik}^{t}} +  \sum_{\fixes{w\in w}}{CM_{wk}^{t} M_{wk}^{t} } + \sum_{e\in E}{H(e) CL_{k}^{t}{O}_{k} X_{ke}^{t}} \right.  \nonumber \\     
&  \left. +  \sum_{e \in E}{E_H(e) CA_{k}^{t} {O}_{k} X_{ke}^{t}} + \sum_{e \in E}{E_F(e) CF_{k}^{t} {O}_{k} X_{ke}^{t}} +  \sum_{s \in S}(CC_{s}^{t}C_{s}^{t}  + CO_{st}O_{s}^{t} )  \nonumber \right\rbrace \\
&  \quad \quad  + \sum_{i \in I, t \in T}  {CH_{i}^{t} I_{i}^{t}} + \sum_{t \in T} CT_{t} \mathcal{E}^t                      - \sum_{k \in K, t \in T}{\alpha_{k}^{t} RF_{k}^{t} } \numberthis \label{FO}
                                 \\
  & \hspace{2mm} \text{s.t:}    \quad (\ref{1.1}),(\ref{1.2}),(\ref{1.3}),(\ref{1.4}),(\ref{1.5}),(\ref{1.6}),(\ref{1.7}),(\ref{1.8}),(\ref{1.9}),(\ref{1.10}),(\ref{1.33}),(\ref{1.11}),(\ref{1.34}),\nonumber \\
  & \hspace{10mm} (\ref{1.12}),(\ref{1.13}),(\ref{1.32}),(\ref{1.14}),(\ref{1.15}),(\ref{1.16}),(\ref{1.17}),(\ref{1.18}),(\ref{1.19}),(\ref{1.20}),\nonumber \\
  & \hspace{10mm} (\ref{1.21}),(\ref{1.22}),(\ref{1.23}),(\ref{1.24}),(\ref{1.25}),(\ref{1.31}),(\ref{1.26}),(\ref{1.27}),(\ref{1.28}),(\ref{1.29}), \nonumber 
\end{align*}
where $E_{H}(e)=\max\{0,$ $H(e)- T(e)\}$ and $E_{F}(e)=\max\{0,$ $T(e)- H(e)\}$ \fixes{are parameters that } denote the \fixes{ number of work shifts that are opened and closed }between period $t-1$ and $t$ \fixes{under state transition $e$}, respectively. The number of workers needed to operate machine $k$ is defined by the parameter $O_k$. \fixes{In \ref{sec:glossary}, a glossary of the used terminology can be found.}

It is easy to show that \textbf{SPT} model without carbon taxes, i.e., the objective function (\ref{FO}) minus \;\;\;  $\sum_{t \in T} CT_{t} \mathcal{E}^t$, is a relaxation of \textbf{SPT} model. We denote this relaxed model as \textbf{SPWT}. Consequently, $Z^*_{SPWT} \leq Z^*_{SPT}$, where $Z^*_{SPWT}$ and $Z^*_{SPT}$ are the optimal solution of \textbf{SPWT} and \textbf{SPT}, respectively. 

\fixes{
\subsection{Effectiveness carbon taxes measures}\label{3.3}

To determine carbon tax effectiveness, we defined two effectiveness measures. 
The first measure is the \textit{transition level}, defined as the proportion of production allocated to each available technology at every period of the planning horizon. 
The second measure is the \textit{transition period}, defined as the period where the weighted proportion of production allocated to clean technologies reaches $\beta$ for the first time, with $\beta \in [0,1]$. 

Let $J$ be the set of available technologies. 
Furthermore, let $N_j$ be an indexed set where $N_j=\{k \in K$ : $k$ is of technology $j\}$, with $\cup_{j \in J} N_j = K$ and $N_j \cap N_p =\emptyset$ for any $j,p \in J: j \neq p$. 
Thus, the transition level of technology $j$ reached in period $t$, is defined as: \begin{align}
    R_j^{t} = \sum_{i \in I}\frac{\displaystyle\sum_{k \in N_j}Y^{t}_{ik}}{\displaystyle\sum_{k \in K}Y^{t}_{ik}} \hspace{10mm} \forall t \in T, j \in J, \label{Rj}
\end{align} 
where $Y^{t}_{ik}$ is an optimal variable of \textbf{SPT} model. 
On the other hand, the transition period where the transition level of clean technologies reaches $\beta$ for the first time, is defined as: 
\begin{align}
     \tau_\beta = \min \{t : \sum_{j \in J} \gamma_j R_j^t  \geq \beta \}, \label{taubeta}
\end{align}
where $\gamma_j$ is the normalized weight for technology $j \in J$, such that $\sum_{j \in J} \gamma_j = 1$ and $\gamma_j > \gamma_p$ if technology $j$ is consider cleaner than technology $p$. We propose a criterion based on the possible emissions as follows: 
\begin{align}
    \gamma_j = \frac{\eta_j^{-1} - \underset{j \in J}{\min}(\eta_j^{-1})}{\displaystyle\sum_{j \in J} \left( \eta_j^{-1} - \underset{j \in J}{\min}(\eta_j^{-1}) \right) } \hspace{10mm} \forall j \in J, 
\end{align}
where, $\eta_j = |N_j|^{-1}\sum_{k \in N_j}\sum_{i \in I} ep_{ik}$ for any $j \in J$. It should be note that the dirtiest technology has a weight equal to $0$. 

\subsection{Evaluation of carbon taxes effectiveness}
A carbon tax is considered useful if a technological transition towards clean production technologies is induced by the end of the planning horizon. 
The effectiveness of the carbon tax is determined by the magnitude of change and the time required to reach $\beta$.   
The weighted average transition level reached by the end of the planning horizon determines the magnitude of the change induced by carbon taxes, i.e. $\sum_{j \in J} \gamma_j R_j^{|T|}$. 
The transition period ($\tau_\beta$) measures how long it takes for the carbon tax to induce the desired technological transition. 
Thus, the higher the transition level and the lower the transition periods, the higher is the effectiveness of the carbon tax in the firm. 
}

\section{Computational study}\label{sec:CS}
In this section, we present our numerical test and results. 
The main objective of the computational study is to explore carbon taxes effectiveness to induce a clean technology transition when the strategic capacity plan considers different technologies, in terms of emissions and investment cost using \textbf{SPT} model. 
Similarly to \cite{drake2016technology} and \cite{song2017capacity}, two types of technologies are considered, more precisely, dirty and clean technology. 

We present an illustrative example of a firm that produces \fixes{eight different natural juice products}, on which carbon taxes are imposed by the environmental authority. 
The strategic capacity planning under carbon taxes considers machines of two technologies, and the \fixes{clean technology transition is illustrated in terms of the transition levels ($R_j^t$). }
\fixes{Moreover, besides the measures proposed in our framework, we compare the total emissions of the \textbf{SPT} model with the \textbf{SPWT} model. }
\fixes{On the other hand, we simulate a firm that faces several exogenous conditions, keeping the endogenous parameters of the industrial example. For this we generate 50 instances, varying the magnitude of the demand. Subsequently, each of these instances were solved considering different relationships between technologies, in terms of emissions and investment costs. }
Thus, to explore carbon taxes effectiveness when the strategic capacity plan considers different technologies, \fixes{the transition level and transition period are computed for each emission and investment relationship. }
\fixes{Each instance considers an eight-year planning horizon. 
However, to represent the continuity of the firm we simulated twelve-years and discard the last four years in the results to }avoid the undesirable effect of the abrupt end of the planning horizon.

All experiments are solved on an Intel Core i5 2.3 GHz x 4 processor with 8 GB RAM, using Gurobi 9.1.1. 
\fixes{We used as termination condition an optimality gap of $10^{-4}$ within a time limit of $36000$ seconds.}

\subsection{Industrial illustrative example}\label{S4.2}
\begin{figure}[h!]
    \centering
    \includegraphics[scale=0.7]{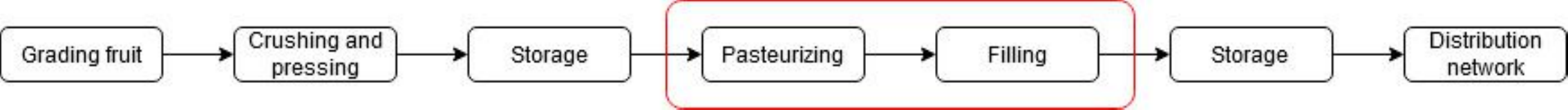}
    \caption{Juice packaging production line.}
    \label{jugo}
\end{figure}

Consider a firm that requires to determine its strategic capacity expansion plan for the upcoming 8 years, during which carbon taxes are implemented. 
The production considers 8 items, which consisting of the juice from 4 different fruits packed in glass containers of $1 000$ and $300$ $cc$. The juice production process consists in the grading of fruit, crushing and pressing the fruit, and storing the juice in stainless-steel containers. Once the container reaches a threshold quantity, the juice is pasteurized at 90 °C and then sent to the filling station. In particular, it is required a capacity expansion for the pasteurizing and filling process, which are considered as a production line. 
Figure \ref{jugo} shows the complete juice-making process\fixes{, where the consider processes are framed.} 

The firm faces a market that grows faster at the beginning of the planning horizon and stabilizes at the end, i.e., it has an increasing demand at a decreasing rate, \fixes{behavior that } represents the life cycle of most products.
Moreover, increasing carbon taxes are imposed throughout the planning horizon. Carbon taxes begin at $\$$35 per carbon emission ton, and increase yearly until reaching to $\$$70 at the end of the planning horizon, \fixes{following } \cite{parry2019}. 

\fixes{For simplicity, let $J=\{1,2\}$, where $j=1$ and $j=2$ are the dirty and clean technology, respectively. }
Let $\widehat{ep}_{ij}$ be the emissions for producing one unit of item $i$ using technology $j$, where $ep_{ik}=\widehat{ep}_{ij}$ for any $i\in I$ and $k\in N_j$.
It is considered that the investment cost ($CI_{k}^{t}$) is higher for the clean technology than for the dirty one. Let $\widehat{CI}_{j}^{t}$ be the investment cost of technology $j$ in period $t$, where $CI_{k}^{t}=\widehat{CI}^{t}_{j}$ for any $k \in N_j$ and $t \in T$.
\fixes{Thus, we consider that carbon emissions and invesment cost are the same for all machines of the same technology. It should be noted that $\gamma_1=0$ and $\gamma_2=1$. }
The parameters used for the industrial example can be found in \ref{param}.

\begin{figure*}[h!]
    \centering 
        \includegraphics[scale=1]{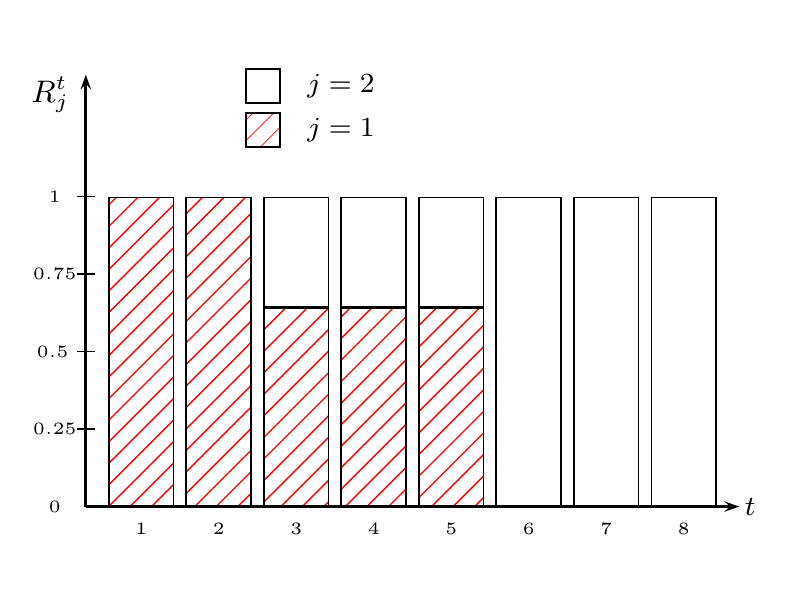}
    \caption{Optimal production allocation by period.}
    \label{fig:production2}
\end{figure*}

\fixes{To illustrate the technological transition induced by carbon taxes, we computed the optimal transition level $R_j^t$ for each $j \in J$ and $t \in T$ }.   
Figure \ref{fig:production2} shows the optimal transition level for both technologies throughout the planning horizon.

From Figure \ref{fig:production2} it can be observed that during the first periods of the planning horizon ($t=1,$ $2$), the firm only uses dirty technology, since the costs of carbon taxes are insufficient for the firm to invest in\fixes{, the more expensive, } clean technology. 
In the following three periods ($t=3,$ $4,$ $5$), the firm starts to use clean technology, but still most of its production is allocated to dirty technology. 
At the beginning of the sixth period ($t=6$), dirty technology \fixes{is no longer used}. Thus, the cost of carbon taxes needed to induce a total technological transition is reached at the beginning of period 6.
Therefore, it can be inferred that carbon taxes become a useful mechanism to induce a technological transition once the carbon emissions costs makes it convenient for the strategic plan to invest in clean machines\fixes{, and the effectiveness not necessarily will be seen with a short-term planning. } 
It should be noted that without carbon taxes all production is allocated to dirty technology, i.e., without carbon taxes there is no incentive for a clean technology transition. 

\begin{figure*}[h!]
        \centering 
        \includegraphics[scale=1]{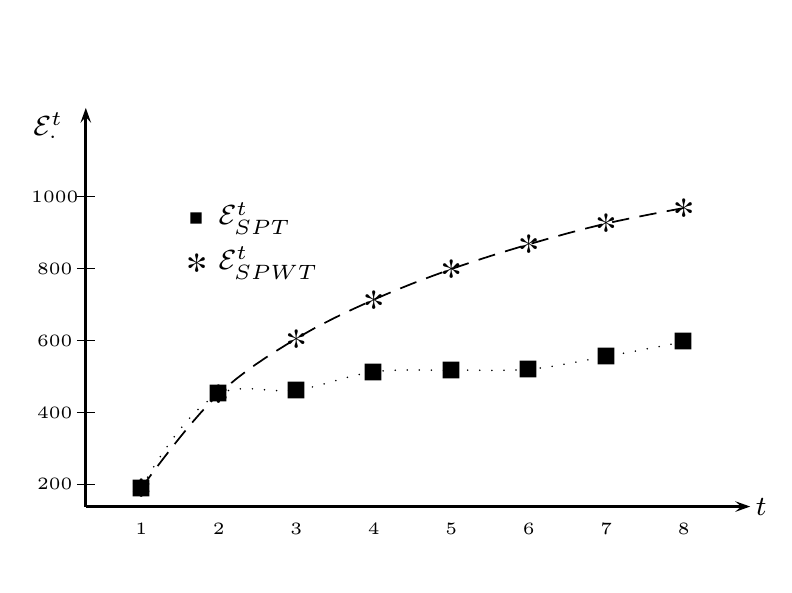}
    \caption{Optimal emissions by period.} 
    \label{fig:emissions2}
\end{figure*}

We computed the total emissions at each period by the optimal solution of \textbf{SPT} and \textbf{SPWT} model respectively.
Figure \ref{fig:emissions2} show the emissions throughout the planning horizon.

From Figure \ref{fig:emissions2} we observe that in the absence of carbon taxes, emissions increase monotonously and grow accordingly with the demand behavior of the firm. On the contrary, under carbon taxes the emissions decouple from the demand. 
Furthermore, it can be observed that once the clean technology is incorporated in production ($t=3$ according to Figure \ref{fig:production2}), the growth rate of emissions decreases considerably. 

\subsection{Test set} 
In what follows, we explore the impact of \fixes{different exogenous parameters, such as the available technology and demand magnitud } on carbon taxes effectiveness to induce a clean technology transition. 
Carbon taxes effectiveness is determined by the clean transition level reached at the end of the planning horizon, i.e., $R_2^8$. 
For the test set, we simplify \textbf{SPT} model, such that the acquisition and replacement of machines become the only variables of interest. \fixes{Thus, } maintenance is removed of \textbf{SPT} model, the work shifts are set to one, and the products are aggregated \fixes{into an equivalent single product}. For each instance of the test set, the aggregated demand of the equivalent product is determine according to $\bar{d_{t}} = \xi \sum_{i \in I} d_{it}$, where $d_{it}$ is the demand used in the illustrative example and $\xi \in \left\lbrace 0.05+\frac{U(0.05,25)-0.05}{10} \right\rbrace$, i.e., $\xi$ varies uniformly between $0.05$ and $2.5$. Similarly, the emissions for the equivalent product are denoted by $\widehat{ep}_{j}$ for any $j \in J$. 
In summary, we study the strategic capacity plan of a single-product firm under carbon taxes that increases yearly from $\$$35 to $\$$70 per carbon emission ton, under different exogenous conditions. 

Two emission relationships ($\frac{\widehat{ep}_{2}}{\widehat{ep}_{1}} \in \{ 0.5, 0.7 \}$), and four investment relationships ($\frac{\widehat{CI}_{2}}{\widehat{CI}_{1}} \in \{ 1.3, 1.4, 1.5, 1.6 \}$) are considered. \fixes{It should be noted that $\frac{\widehat{ep}_{2}}{\widehat{ep}_{1}} < 1$ and $\frac{\widehat{CI}_{2}}{\widehat{CI}_{1}} > 1$, since we assume that the less emitter technology is more expensive. } 

Therefore, each instance of the test set is evaluated for the eight combinations of emissions and investment relationships\fixes{, generating $400$ scenarios. }The parameters used for the test set are the aggregated parameters of the illustrative example. 

\subsubsection{Transition level at the end of the planning horizon}

\begin{figure}[h!]
    \vspace{3.5mm}
    \centering
    \begin{subfigure}[b]{1\textwidth}
    \centering 
    \includegraphics[scale=0.7]{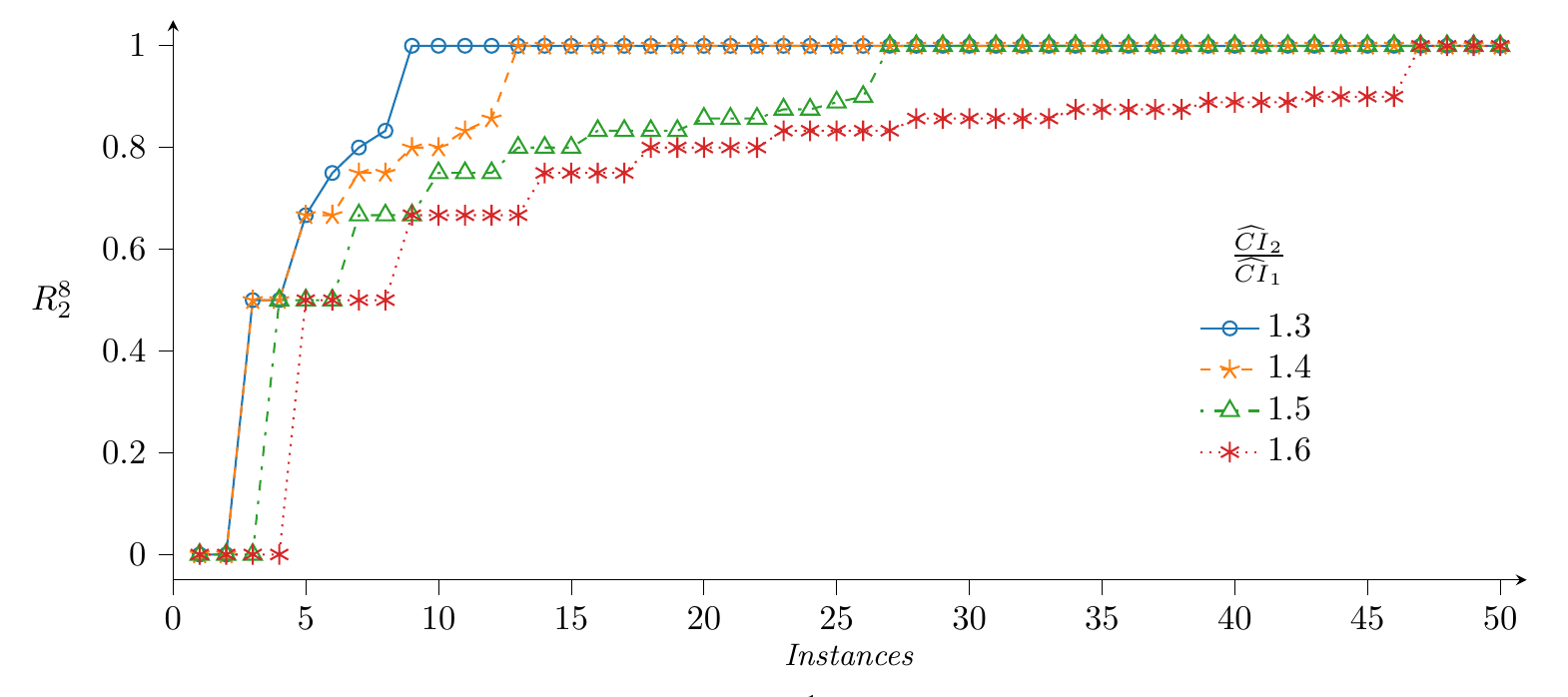}
    \caption{$\frac{\widehat{ep}_2}{\widehat{ep}_1}=0.5$.}
    \label{plot05}
    \end{subfigure}

    \centering
    \vspace{3.5mm}
    \begin{subfigure}[b]{1\textwidth}
    \centering
    \includegraphics[scale=0.7]{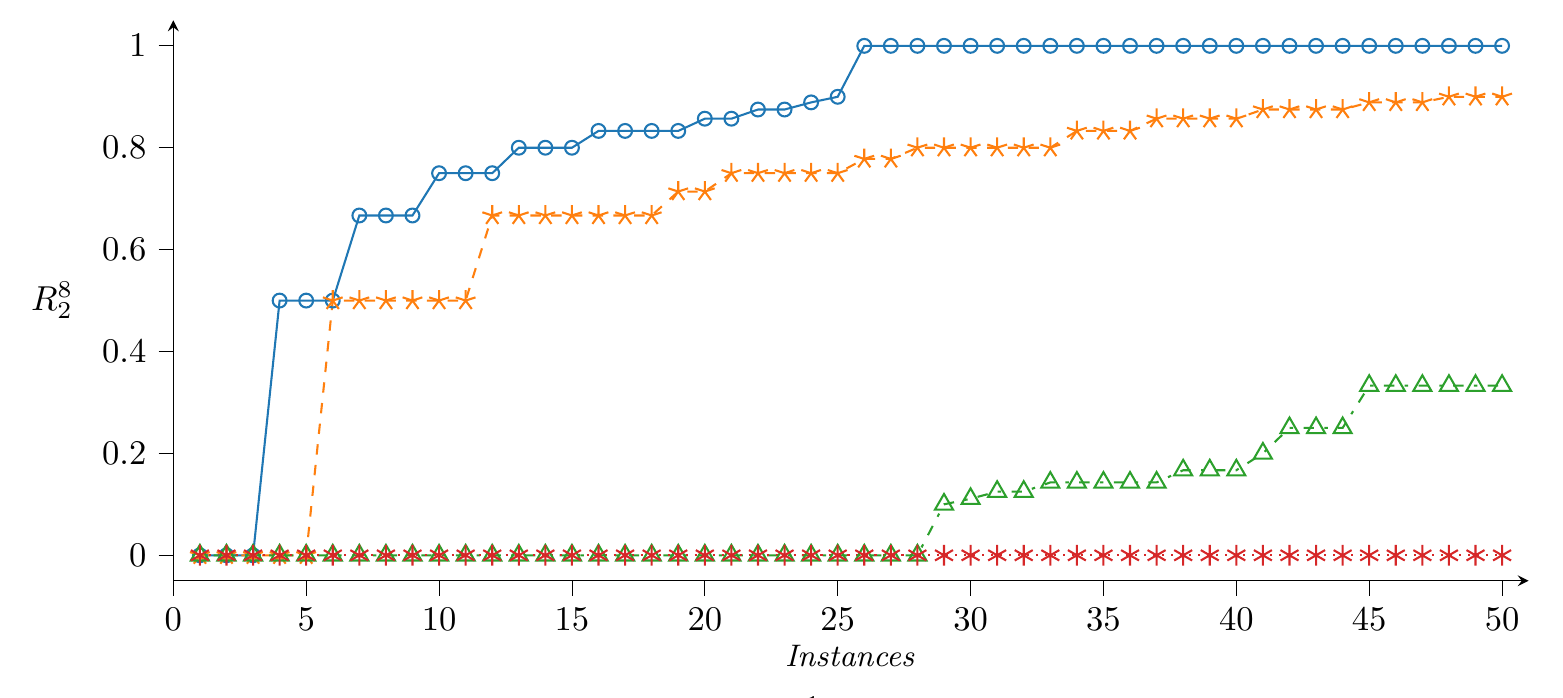} 
    \caption{$\frac{\widehat{ep}_2}{\widehat{ep}_1}=0.7$.}
    \label{plot07}
    \end{subfigure}
    \caption{$R_2^8$ for different emissions and investment relationships.}\label{acumulado}
\end{figure}

To study the effect of the exogenous conditions considered, we analyze the clean transition level reached at the end of the planning horizon. 
Figure \ref{acumulado} shows the proportion of allocated production to clean technology at the end of the planning horizon ($R_2^8$) \fixes{for every scenario.} 

From Figure \ref{acumulado} we observe that carbon taxes effectiveness is highly influenced by the emission relationship between the clean and dirty technology\fixes{, and the demand}.
More precisely, the greater the emission reduction of the clean technology, the greater the probability of reaching a technological transition. In addition, the higher the investment cost of the clean technology, the lower the probability of reaching a technological transition. \fixes{Regarding the demand, the lower the demand size, the lower is the probability of reaching a technological transition under any technology relationship. }
Figure \ref{plot05} shows that the probability of achieving a total technological transition is greater than zero for every investment relationship considered.
On the other hand, Figure \ref{plot07} shows that the probability of achieving a total technological transition is greater than zero only for $\frac{\widehat{CI}_2}{\widehat{CI}_1}=1.3$, and the probability of achieving a technological transition is even zero for $\frac{\widehat{CI}_2}{\widehat{CI}_1}=1.6$. 
Table \ref{probab} summarizes the probabilities for achieving different levels of technological transition, and the expected value of the transition at the end of the planning horizon for each emission and investment relationship.

\begin{table}[htbp]
  \centering
  \resizebox{\textwidth}{!}{  
    \begin{tabular}{cccccccccccc}
    \toprule
          & \multicolumn{5}{c}{$\frac{\widehat{ep}_2}{\widehat{ep}_1}=0.5$}         &       & \multicolumn{5}{c}{$\frac{\widehat{ep}_2}{\widehat{ep}_1}=0.7$} \\
\cmidrule{2-6}\cmidrule{8-12}    $\frac{\widehat{CI}_2}{\widehat{CI}_1}$ & $\mathbb{P}(R_2^8=0)$ & $\mathbb{P}(R_2^8 \geq 0.50)$ & $\mathbb{P}(R_2^8\geq 0.75)$ & $\mathbb{P}(R_2^8\geq 1)$ & $\mathbb{E}(R_2^8)$ &   & $\mathbb{P}(R_2^8=0)$ & $\mathbb{P}(R_2^8 \geq 0.50)$ & $\mathbb{P}(R_2^8\geq 0.75)$ & $\mathbb{P}(R_2^8\geq 1)$ & $\mathbb{E}(R_2^8)$ \\
    \midrule
    1.3   & 0.04  & 0.96  & 0.90  & 0.82  & 0.92  &       & 0.06  & 0.94  & 0.82  & 0.56  & 0.85 \\
    1.4   & 0.06  & 0.94  & 0.86  & 0.76  & 0.89  &       & 0.10   & 0.90  & 0.62  & 0.06  & 0.70 \\
    1.5   & 0.06  & 0.94  & 0.82  & 0.52  & 0.84  &       & 0.56  & 0.00  & 0.00  & 0.00  & 0.22 \\
    1.6   & 0.08  & 0.92  & 0.74  & 0.08  & 0.74  &       & 1.00     & 0.00  & 0.00  & 0.00  & 0.00 \\
    \bottomrule
    \end{tabular}
    }
  \caption{Transition expectancy values.} \label{probab}
\end{table}

Table \ref{probab} shows that for $\frac{\widehat{ep}_2}{\widehat{ep}_1}=0.5$, the probability to achieve a transition is strictly higher than for $\frac{\widehat{ep}_2}{\widehat{ep}_1}=0.7$. 
Moreover, for $\frac{\widehat{ep}_2}{\widehat{ep}_1}=0.5$ and $\frac{\widehat{ep}_2}{\widehat{ep}_1}=0.7$, the expected values are $\mathbb{E}(R^8_2)=0.74$ and $\mathbb{E}(R^8_2)=0$ respectively.
Hence, we infer that the emission relationship has a great impact on carbon taxes effectiveness to induce a technological transition. 
Furthermore, when $\frac{\widehat{ep}_2}{\widehat{ep}_1}=0.5$ the investment relationship has a lower impact, e.g., $\mathbb{P}(R_2^8 \geq 0.50) = 0.96$ for  $\frac{\widehat{CI}_2}{\widehat{CI}_1} = 1.3$ and $\mathbb{P}(R_2^8 \geq 0.50) = 0.92$ for  $\frac{\widehat{CI}_2}{\widehat{CI}_1} = 1.6$.
On the contrary, when  $\frac{\widehat{ep}_2}{\widehat{ep}_1}=0.7$ the investment relationship has a higher impact on the transition, e.g., $\mathbb{P}(R_2^8 \geq 0.50) = 0.94$ for  $\frac{\widehat{CI}_2}{\widehat{CI}_1} = 1.3$ and $\mathbb{P}(R_2^8 \geq 0.50) = 0$ for  $\frac{\widehat{CI}_2}{\widehat{CI}_1} = 1.6$.
Hence, we infer that the investment relationship has a greater impact on carbon taxes effectiveness to induce a technological transition when the emission relationship is higher. 
\fixes{In summary, when clean technology is more expensive and does not reduce emissions enough and the demand is small, the firm has low or none incentive to change from a dirty technology, since it is more convenient to pay the carbon tax than to invest in a clean technology.}

\subsubsection{Effect of technology relations on transition periods}

\begin{figure}[h!]
    \centering
    \begin{subfigure}[b]{0.4\textwidth}
    \centering
    \includegraphics[scale=0.71]{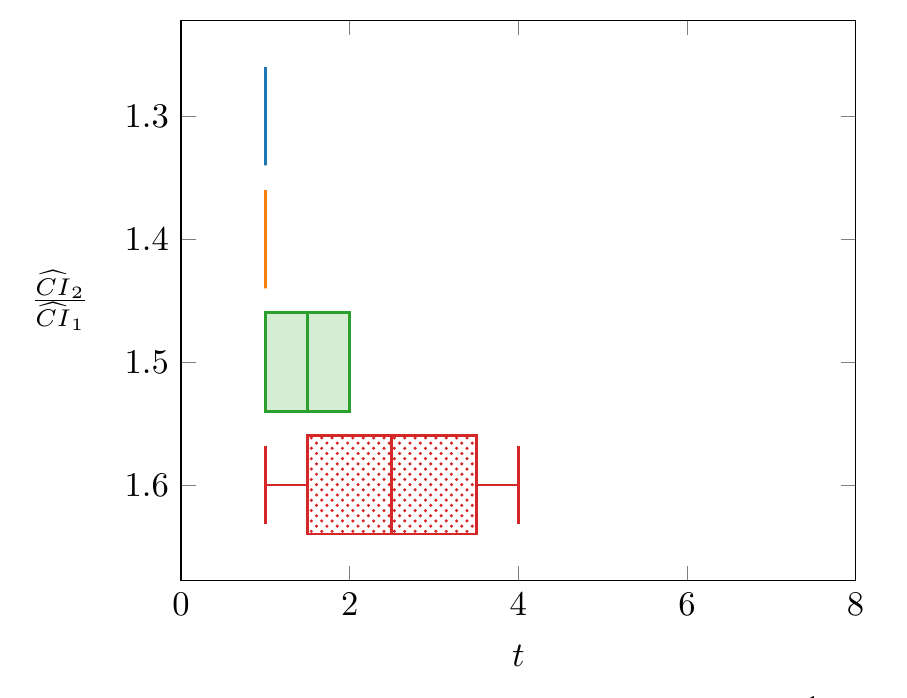}
    \caption{$\frac{\widehat{ep}_2}{\widehat{ep}_1}=0.5$,$\; \tau_{0.50}$.}
    \label{box1}
    \end{subfigure}
    \hspace{5mm}
    \begin{subfigure}[b]{0.4\textwidth}
    \centering
    \includegraphics[scale=0.71]{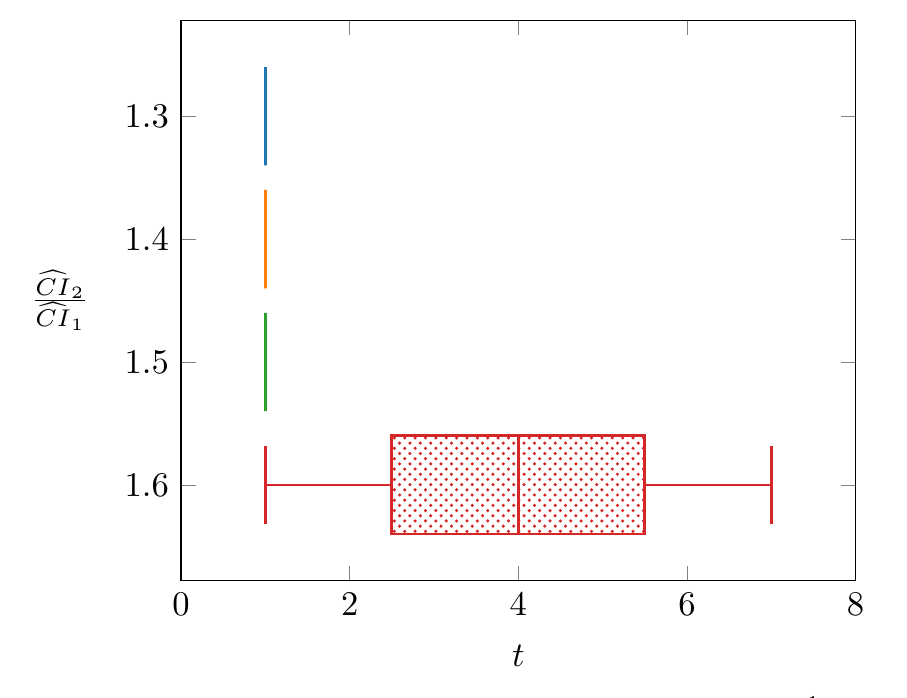}
    \caption{$\frac{\widehat{ep}_2}{\widehat{ep}_1}=0.5$,$\; \tau_{0.75}$.}
    \label{box2}
    \end{subfigure}

    \centering
    \vspace{5mm}
    \begin{subfigure}[b]{0.4\textwidth}
    \centering
    \includegraphics[scale=0.71]{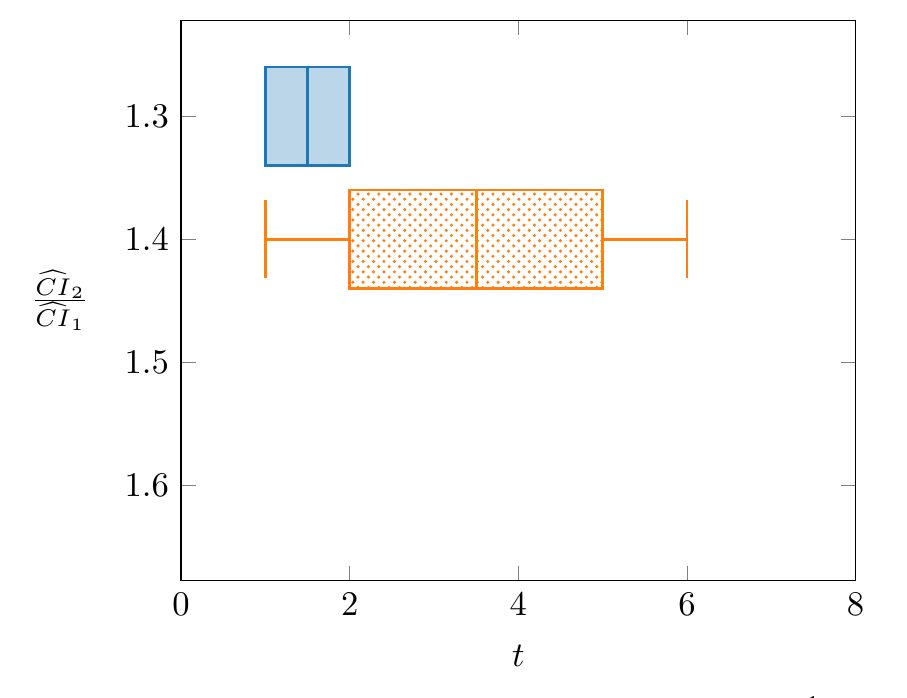}
    \caption{$\frac{\widehat{ep}_2}{\widehat{ep}_1}=0.7$,$\; \tau_{0.50}$.}
    \label{box3}
    \end{subfigure}
    \hspace{5mm}
    \vspace{5mm}
    \begin{subfigure}[b]{0.4\textwidth}
    \centering
    \includegraphics[scale=0.71]{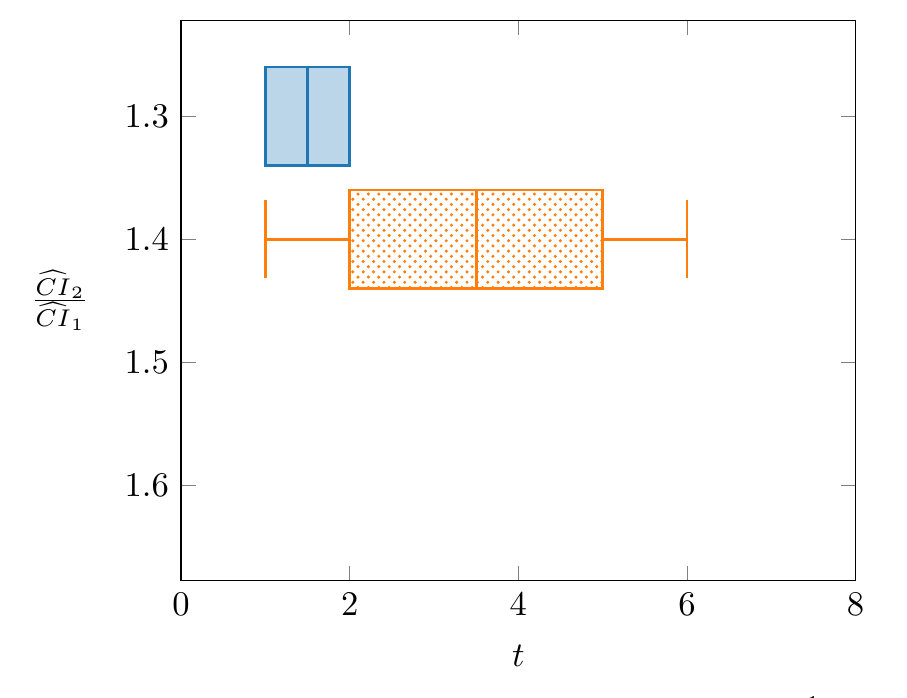}
    \caption{$\frac{\widehat{ep}_2}{\widehat{ep}_1}=0.7$,$\; \tau_{0.75}$.}
    \label{box4}
    \end{subfigure}
    \caption{$\tau_\beta$ for different emissions and investment relationships.}\label{boxes}
\end{figure}

\fixes{For every scenario analized, } we computed $\tau_\beta$, for $\beta \in \{0.50, 0.75\}$.
For each $\beta$, only scenarios where $\tau_\beta < \infty$ are considered, e.g., for $\frac{\widehat{ep}_2}{\widehat{ep}_1}=0.5$ and $\frac{\widehat{CI}_2}{\widehat{CI}_1}=1.4$ the $94\%$ of instances reach or surpass $\beta=0.50$, according to Table \ref{probab}. 
Figure \ref{boxes} shows the box-plots of $\tau_\beta$ for $\beta \in \{0.50,0.75\}$ and the different emissions and investment relationships. 

From Figure \ref{box1}, we observe that for $\frac{\widehat{CI}_2}{\widehat{CI}_1} \in \{1.3,1.4\}$, all instances that achieved $\beta=0.50$ reach it at the beginning of the planning horizon, i.e., $\tau_{0.50}=1$.
In contrast, for $\frac{\widehat{CI}_2}{\widehat{CI}_1}=1.5$ and $1.6$, we observe that $\tau_{0.50}=1.5$ and $\tau_{0.50}=2.5$, respectively. 
From Figure \ref{box2}, we observe that for $\frac{\widehat{CI}_2}{\widehat{CI}_1} \in \{1.3,1.4,1.5\}$, all instances which achieved $\beta=0.75$ reach it at $\tau_{0.75}=1$.
In contrast, for $\frac{\widehat{CI}_2}{\widehat{CI}_1}=1.6$, we observe that $\tau_{0.75}=4$. 
On the other hand, Figure \ref{box3} only contains the box-plots of $\frac{\widehat{CI}_2}{\widehat{CI}_1} \in \{1.3,1.4\}$, where $\tau_{0.50}=1.5$ and $\tau_{0.50}=3.5$, respectively. 
Further, Figure \ref{box4} only contains the box-plots of $\frac{\widehat{CI}_2}{\widehat{CI}_1} = 1.3$ and $1.4$, where $\tau_{0.75}=1.5$ and $\tau_{0.75}=4$, respectively.
It should be noted that $\frac{\widehat{CI}_2}{\widehat{CI}_1}\in \{1.5,1.6\}$, no instance reaches $\beta \in \{ 0.5, 0.7 \}$. Consequently, there is no box-plot for $\frac{\widehat{CI}_2}{\widehat{CI}_1}\in \{1.5,1.6\}$ under $\frac{\widehat{ep}_2}{\widehat{ep}_1}=0.7$.
\fixes{In summary, we observe that the firm has incentive to transition earlier to a cleaner technology when the benefit from reducing emissions is high and the increase of the investment cost is low, i.e., $\tau_{\beta}$ decreases when the emission and investment relationships decrease.}

\section{Conclusion}\label{sec:CO}

\fixes{
In this paper, we introduce a framework to evaluate carbon taxes effectiveness to induce a clean technology transition on individual firms. 
The framework is composed by three main elements. 
The first element is a strategic capacity planning under carbon taxes model, which determines the optimal strategic plan in terms of acquisition and usage of technologies throughout the planning horizon. Thus, determines whether there is a technology transition, and if so, the optimal technological transition strategy.  
The second element correspond to the measures to evaluate the carbon tax. We define the transition level and the transition period, as the proportion of allocated production to each technology and the period where the desired clean transition level is reached, respectively.  
The third element correspond to the evaluation of the carbon tax effectiveness itself, using the defined measures. } 

\fixes{We apply the framework to an industrial illustrative example where the strategic capacity plan under carbon taxes faces a clean and a dirty technology. 
The technological transition is illustrated in terms of the transition levels for each technology and total emissions throughout the planning horizon. }
Several instances, based on the industrial illustrative example, were generated to explore the impact of technology relationships, in terms of emissions and investment, and size of the demand on carbon taxes effectiveness. From our computational experiments, we observe that a full clean technology transition is not reached in the planning horizon when the emissions of the clean technology are a 70$\%$ of the dirty technology, and the investment cost of clean technologies is 140$\%$ or higher than the dirty technology. 
Thus, we infer that for costly clean technology that does not greatly reduce emissions (compared to dirty technology), \textcolor{black}{carbon taxes are not effective on incentivizing the firm to change technologies. Furthermore, we infer that a firm that face small demand has lower incentive to acquire a more costly technology. } 
Moreover, carbon taxes alone do not ensure a clean technology transition because their effectiveness is highly influenced by the incentive of firms to invest in a more expensive technology to reduce emissions.  

\textcolor{black}{There are several issues left for future research. 
Since carbon taxes effectiveness is affected by the demand, i.e., the size of the firm, the first issue to be explored is to determine the optimal carbon tax policy for each firm size, such that all firms in the industry are encouraged to a clean technology transition.
The environmental authority can choose from a range of mechanism to reduce emissions and incentivize a transition to clean technologies. Hence, the second issue is to incorporate into the framework other common mechanism, such as the emissions trading system and subsidies, allowing the comparison of the effectiveness of different mechanisms in inducing a clean technology transition in a firm. 
It is reasonable to consider that demand and available technology is unknown for the firm throughout the planning horizon, due to changes in the market and the incentive promoted by carbon taxes for suppliers to lower technologies emissions \citep{boyce2018carbon}. 
Following this, the third issue is to consider a stochastic strategic capacity planning model in the framework, such that the model consider uncertainty in the demand and available technology. 
}
The last issue left is the scalability of \textbf{SPT} model, as it becomes considerably harder to prove the optimality of the solution with a real case complexity. 


\section*{Acknowledgements}
\textcolor{black}{
Professor Pablo Escalona is grateful for the support of ANID, through grant FONDECYT 11200287. Jorge Weston has received funding from the European Union’s Horizon 2020 research and innovation program under the Marie Skłodowska-Curie grant agreement No 764759 ETN “MINOA”.  
}

\bibliographystyle{elsarticle-harv}\biboptions{authoryear}
\bibliography{References}

\newpage
\appendix
\allowdisplaybreaks

\section{Glossary of terms} \label{sec:glossary}
\begin{table}[!htbp]
\centering{
\scalebox{0.70}{
\begin{tabular}{c l}  \hline                             
Sets & Definition  \\ \hline
$T$ & Set of  periods indexed by $t$, with $t=1,..., |T|$ \\
$I$ & Set of items indexed by $i$, with $i=1,...,|I|$ \\
$K^0$ & Set of operative machines at the beginning of the planning horizon  \\
$K^c$ & Set of machine that can be acquired over the planning horizon \\
$K$ & Set of machine indexed by $k=1,...,N$, $K=K^0 \cup K^c$ \\
$S$ & Set of states in which machine can be situated over the planning horizon, $S=\{0,1,2,3\}$ where \\
& $s=0$ the machine is inoperative; $s=1$ the machine is operative to be used one work shift; \\
& $s=2$ the machine is available to be used two work shifts; $s=3$ the machine is available to be used three work shifts; \\
$E$ & Set of state transitions, with $E\varsubsetneq \{ (s,s\textprime) \in S \times S\}$ \\ 
\hline
Operational parameters &  \\ 
\hline
$d_{it}$ & Demand for item $i$ in period $t$\fixes{, measure in units.}\\ 
$s_0$ & State of machine $k$ at the beginning of the planning horizon ($t=0$), for any $k\in K^0$  \\
$H_i^{0}$ & On-hand inventory of item $i$ at the beginning of period $1$\fixes{, measure in units.} \\
$FTM_k$ & Fixed time maintenance of machine $k$\fixes{, measure in units of time.} \\ 
$RMT_{wk}$ & Time required for the $w$\textit{-th} maintenance of machine $k$\fixes{, measure in units of time.}  \\
$r_{ik}$ & Production rate of item $i$, produced with machine $k$\fixes{, measure in unit per unit of time.} \\
$\mu_{k}$ & Maximum utilization for each machine $k$ \\
$v_k$ & Useful life of machine $k$\fixes{, measure in units of time.}\\ 
$l$ & Available working time for any work shift $s$\fixes{, measure in units of time.}\\ 
$ep_{ik}$  & Emissions for producing one unit of item $i$ using machine $k$\fixes{, measure in $co_2$ ton per unit.} \\ 
$eh_{i}$ & Emissions of holding per unit per unit of item $i$\fixes{, measure $co_2$ ton per unit of unit time.}\\
$O_k$ & Number of workers needed to operate machine $k$  \\ 
\hline
Cost parameters  &  \\ 
\hline
$CI_{k}^{t}$ & Investment cost of acquiring a machine $k$ in period $t$\fixes{, measure in monetary units.}\\ 
$CP_{ik}^{t}$ & Production cost for one item $i$  produced by machine $k$ in the period $t$\fixes{, measure in monetary units per unit.} \\ 
$CM_{wk}^{t}$ & Preventive $w$-th maintenance cost of machine $k$ in period $t$\fixes{, measure in monetary units per maintenance.} \\ 
$CL_k^t$ & Labor cost for machine $k$ in period $t$\fixes{, measure in monetary units per worker.} \\ 
$CA_{k}^{t}$ & Cost of hiring a worker for machine $k$ at the beginning of period $t$ \fixes{, measure in monetary units per worker.}\\
$CF_{k}^{t}$ &  Cost of firing a worker of machine $k$ at the beginning of period $t$, monetary unit per worker\fixes{, measure in monetary units per worker.}\\
$CO_{s}^{t}$ & Cost of opening work shift $s$ in period $t$\fixes{, measure in monetary units.} \\
$CC_{s}^{t}$ & Cost of closing work shift $s$ in period $t$\fixes{, measure in monetary units.}\\
$CH_{i}^{t}$ & Holding cost per unit and unit time of item $i$ in period $t$\fixes{, measure in monetary units per unit and unit time.} \\
$CT_{t}$ & Value of the carbon tax in period $t$\fixes{, measure in monetary units per $co_2$ ton.}\\
$\alpha_{k}^{t}$ & Price per hour of useful life remaining of used machine $k$ in period $t$\fixes{, measure in monetary units per unit time.} \\
\hline                       
Parameter functions & \\ \hline
$H(e)$ & State of machine $k$ in period $t$ \\
$T(e)$ & State of machine $k$ in period $t-1$ \\ 
$E_{H}(e)$ & Increase in work shifts in period $t$, $E_{H}(e)=\max\{0,H(e)- T(e)\}$ \\
$E_{F}(e)$ & Decrease in work shifts in period $t$, $E_{F}(e)=\max\{0,T(e)-H(e)\}$ \\ 
$\mathcal{E}^t$ & Total emissions in period $t$, $\mathcal{E}^t = \sum_{i \in I}(\sum_{k \in K} {ep_{ik}  Y_{ik}^{t}} + eh_{i} I_{i}^{t})$\fixes{, measure in $co_2$ ton.}\\
\hline

Variables &  \\ \hline
$X_{ke}^{t}$ & $1$ if transition state $e$ occurs for machine $k$ at the beginning of period $t$; $0$ otherwise \\
$Z_{s}^{t}$ & $1$ if $s$ work shifts are utilized in period $t$; $0$ otherwise \\
$O_{s}^{t}$ & $1$ if work shift $s$ is opened in period $t$; $0$ otherwise\\
$C_{s}^{t}$ & $1$ if work shift $s$ is closed in period $t$; $0$ otherwise \\
$Y_{ik}^{t}$ & Units of item $i$ produce with machine $k$ in period $t$\fixes{, measure in units.} \\
$I_{i}^{t}$ & Units of item $i$ in inventory at the end of period $t$\fixes{, measure in units.}\\
$M_{wk}^{t}$ & $1$ if maintenance $w$-th is done to machine $k$ in period $t$; $0$ otherwise \\
$TM_{k}^{t}$ & Accumulated production time of machine $k$ from its last maintenance at the beginning of period $t$\fixes{, measure in units of time.}\\
$RL_{k}^{t}$ & Remaining useful life of machine $k$ at the beginning of period $t$\fixes{, measure in units of time.}\\
$RF_{k}^{t}$ & Residual life of machine $k$ when sold in period $t$\fixes{, measure in units of time.} \\
 \hline
 
Subset partitions &  \\ \hline
$E^0$ &  Transition that represent a machine that remain inoperative, $E^0 = \{e_{0} \}$ \\
$E^1$ & Transitions that represent the entry into operation of a machine, $E^1 = \{e_{1},e_{2},e_{3} \}$ \\
$E^2$ & Transitions that represent a machine that was and stays in operation, $E^2=\{e_5,e_6,e_7,e_{9},e_{10},e_{11},e_{13},e_{14},e_{15} \}$ \\
$E^3$ &
Transitions that represent the discard of a machine, $E^3 = \{e_4,e_8,e_{12}\}$ \\
\hline 

\end{tabular}}}
\caption{Glossary of terms.}
\label{tab: tabla1}
\end{table}

\newpage
\section{Parameters of the illustrative example} \label{param}

\newcommand{\localtextbulletone}{\textcolor{black}{\raisebox{.45ex}{\rule{.6ex}{.6ex}}}}
\renewcommand{\labelitemi}{\localtextbulletone}

\begin{itemize}
	\item \textit{Production rate}: $\;\; r_{ik}=r_{i} \;\; \forall k \in K$
	\FloatBarrier
    \begin{table}[h!]
  \centering
  \resizebox{.5\linewidth}{!}{
    \begin{tabular}{lcccccccc}
    \toprule
          & \textit{i=1}   & \textit{i=2}   & \textit{i=3}   & \textit{i=4}   & \textit{i=5}   & \textit{i=6}   & \textit{i=7}   & \textit{i=8} \\
    \midrule
    $r_i$    & 480   & 672   & 576   & 336   & 528   & 624   & 768   & 816 \\
    \bottomrule
    \end{tabular}%
    }
\end{table}%

	\item \textit{Utilization}: ${\mu_k}=0.85 \; \forall k \in K$
	
	\item \textit{Useful life}: $v_k=20000$ $\forall k \in K$
	
	\item \textit{Time by shift}: $l=2080$
	
	\item \textit{Fixed time maintenance}: $FTM_k=5000$ $\forall k \in K$
	
	\item \textit{Time by maintenance}: $RMT_{wk}=4$ $\forall w \in W, k \in K$  
	
	\item \textit{Workers by machine}: $O_k=2 \; \forall k \in K$
	
	\item \textit{Shift work at $t=0$}: $s_0=0$
	
	\item \textit{Inventory emissions}: 
	\FloatBarrier
    \begin{table}[h!]
    \centering
    \hspace{7mm}
      \resizebox{.55\linewidth}{!}{
    \begin{tabular}{lrrrrrrrr}
    \toprule
          & \multicolumn{1}{c}{$i=1$} & \multicolumn{1}{c}{$i=2$} & \multicolumn{1}{c}{$i=3$} & \multicolumn{1}{c}{$i=4$} & \multicolumn{1}{c}{$i=5$} & \multicolumn{1}{c}{$i=6$} & \multicolumn{1}{c}{$i=7$} & \multicolumn{1}{c}{$i=8$} \\
    \midrule
    $eh_i$    & 0.023 & 0.032 & 0.027 & 0.016 & 0.025 & 0.029 & 0.036 & 0.038 \\
    \bottomrule
    \end{tabular}%
    }
    \end{table}
	
	\item \textit{Production cost}: $\;\; CP_{ik}^t=CP_{i}^t \;\; \forall k \in K$
	\FloatBarrier
    \begin{table}[h!]
    \centering
    \hspace{7mm}
    \resizebox{.6\linewidth}{!}{
    \begin{tabular}{lrrrrrrrr}
    \toprule
          & \multicolumn{1}{c}{$i=1$} & \multicolumn{1}{c}{$i=2$} & \multicolumn{1}{c}{$i=3$} & \multicolumn{1}{c}{$i=4$} & \multicolumn{1}{c}{$i=5$} & \multicolumn{1}{c}{$i=6$} & \multicolumn{1}{c}{$i=7$} & \multicolumn{1}{c}{$i=8$} \\
    \midrule
    $CP^1_i$    & 0.075 & 0.105 & 0.09  & 0.053 & 0.083 & 0.098 & 0.12  & 0.128 \\
    \bottomrule
    \end{tabular}
    }
    \end{table}
    
	\item \textit{Labor cost}: $CL_k^1=4500$ $\forall k \in K$
	
	\item \textit{Cost of hiring and firing}: $CA_k^1=CF_k^1=5000$ $\forall k \in K$
	
	\item \textit{Cost of opening and closing a shift}: $CT_s^1=CU_s^1=5000$ $\forall s \in S$
	
	\item \textit{Holding cost}:
	\FloatBarrier
    \begin{table}[h!]
  \centering
      \resizebox{.55\linewidth}{!}{
    \begin{tabular}{lrrrrrrrr}
    \toprule
          & \multicolumn{1}{c}{$i=1$} & \multicolumn{1}{c}{$i=2$} & \multicolumn{1}{c}{$i=3$} & \multicolumn{1}{c}{$i=4$} & \multicolumn{1}{c}{$i=5$} & \multicolumn{1}{c}{$i=6$} & \multicolumn{1}{c}{$i=7$} & \multicolumn{1}{c}{$i=8$} \\
    \midrule
    $CH_i^1$    & 0.66  & 0.92  & 0.79  & 0.46  & 0.72  & 0.85  & 1.05  & 1.12 \\
    \bottomrule
    \end{tabular}%
    }
\end{table}%
	
	\item \textit{Selling factor}: $\alpha_k^1=\dfrac{0.8 CI_k^1}{v_k}$ $\forall k \in K$
	
	\item \textit{Discount factor}: The discount factor used is $10\%$. Every costs behave over time as $C_{t}=C_{1}(1+0.1)^{1-t}$
	
	\item \textit{Emissions of production}: 
	
	\FloatBarrier
    \begin{table}[h!]
    \centering
    \hspace{8mm}
    \resizebox{.55\linewidth}{!}{
    \begin{tabular}{lcccccccc}
    \toprule
          & $i=1$   & $i=2$   & $i=3$   & $i=4$   & $i=5$   & $i=6$   & $i=7$   & $i=8$ \\
    \midrule
    $\widehat{ep}_{i1}$   & 0.30  & 0.21  & 0.25  & 0.42  & 0.27  & 0.23  & 0.18  & 0.17 \\
    $\widehat{ep}_{i2}$   & 0.30  & 0.21  & 0.25  & 0.42  & 0.27  & 0.23  & 0.18  & 0.17 \\
    \bottomrule
    \end{tabular}%

    }
\end{table}
	
	\item \textit{Investment cost}: $\widehat{CI}_1^1=65000 \;$ and $\widehat{CI}_2^1=104000 $
	
	\item \textit{Maintenance cost}: $\widehat{CM}_{wk}^1=600 \;$ 
	
    \item \textit{Demand:}
    \FloatBarrier
    \begin{table}[h]
        \centering
        \resizebox{.95\linewidth}{!}{
        
    \begin{tabular}{ccccccccccccc}
    \toprule
    $d_{it}$ & $t=1$   & $t=2$   & $t=3$   & $t=4$   & $t=5$   & $t=6$   & $t=7$   & $t=8$   & $t=9$   & $t=10$  & $t=11$  & $t=12$ \\
    \midrule
    $i=1$   & 20000 & 47726 & 63945 & 75452 & 84378 & 91670 & 97836 & 103178 & 107889 & 112103 & 115916 & 119396 \\
    $i=2$   & 28000 & 66816 & 89522 & 105633 & 118129 & 128339 & 136971 & 144449 & 151045 & 156945 & 162282 & 167155 \\
    $i=3$   & 24000 & 57271 & 76733 & 90542 & 101253 & 110004 & 117404 & 123813 & 129467 & 134524 & 139099 & 143276 \\
    $i=4$   & 14000 & 33408 & 44761 & 52816 & 59064 & 64169 & 68485 & 72224 & 75522 & 78472 & 81141 & 83577 \\
    $i=5$   & 22000 & 52498 & 70339 & 82997 & 92815 & 100837 & 107620 & 113495 & 118678 & 123314 & 127507 & 131336 \\
    $i=6$   & 26000 & 62044 & 83128 & 98087 & 109691 & 119172 & 127187 & 134131 & 140256 & 145734 & 150691 & 155215 \\
    $i=7$   & 32000 & 76361 & 102311 & 120723 & 135004 & 146673 & 156538 & 165084 & 172622 & 179365 & 185465 & 191034 \\
    $i=8$   & 34000 & 81134 & 108706 & 128268 & 143442 & 155840 & 166322 & 175402 & 183411 & 190576 & 197057 & 202974 \\
    \bottomrule
    \end{tabular}%
        }
        \end{table}
    
\end{itemize}

\end{document}